\documentclass[12pt,a4paper]{article} 

\usepackage{a4}

\usepackage[intlimits,tbtags]{amsmath} 
\usepackage{amsfonts,amssymb,amsthm}
\usepackage{amscd}

\usepackage{cite,latexsym} 

\newcommand{\E}{\mathrm{e}}

\newcommand{\tr}{\triangleright}
 
\newcommand{\R}{\mathcal{R}}
\newcommand{\F}{\mathcal{F}}
\newcommand{\Fs}{{\mathcal{F}_{\!\mathrm{s}}}}

\newcommand{\id}{\mathrm{id}}

\newcommand{\Xcal}{\mathcal{X}}

\newcommand{\Ocal}{\mathcal{O}}

\newcommand{\Tcal}{\mathcal{T}}

\newcommand{\Ket}[1]{\lvert #1\rangle} 
\newcommand{\Bra}[1]{\langle #1\rvert}

\newcommand{\lrAngle}[1]{\langle #1\rangle}

\newtheorem{Definition}{Definition} 
\newtheorem{Theorem}{Theorem}
\newtheorem{Proposition}{Proposition}
\newtheorem{Corollary}{Corollary}

\newcommand{\h}{\hbar}

\newcommand{\Usu}{{\mathcal{U}(\mathrm{su}_2)}}

\newcommand{\Usl}{{\mathcal{U}(\mathrm{sl}_2)}}

\newcommand{\Uhsu}{{\mathcal{U}_\h(\mathrm{su}_2)}}
\newcommand{\Uhg}{{\mathcal{U}_\h(\mathfrak{g})}}
\newcommand{\Uqg}{{\mathcal{U}_q(\mathfrak{g})}}
\newcommand{\Ug}{{\mathcal{U}(\mathfrak{g})}}

\newcommand{\UhslC}{{\mathcal{U}_{\h}(\mathrm{sl}_2(\mathbb{C}))}}
\newcommand{\UslC}{{\mathcal{U}(\mathrm{sl}_2(\mathbb{C}))}}
\newcommand{\slC}{{\mathrm{sl}_2(\mathbb{C})}}
\newcommand{\so}{{\mathrm{so}_4}}
\newcommand{\Uhso}{{\mathcal{U}_{\h}(\mathrm{so}_4)}}
\newcommand{\Uso}{{\mathcal{U}(\mathrm{so}_4)}}

\newcommand{\SUh}{{SU_{\h}(2)}} 

\newcommand{\Mink}{{\Xcal(\mathbb{R}^{1,3})}}
\newcommand{\Minkh}{{\Xcal_\h(\mathbb{R}^{1,3})}}

\newcommand{\CGs}[6]{
  \bigl(\begin{smallmatrix}#1\! &  #2\vphantom{#3} \\
    #4\! & #5\vphantom{#6} \end{smallmatrix} \!\!\bigm|\!\!
  \begin{smallmatrix}
    \vphantom{#1#2}#3 \\ \vphantom{#4#5}#6
  \end{smallmatrix}\bigr)}

\newcommand{\CGqs}[6]{\CGs{#1}{#2}{#3}{#4}{#5}{#6}_{\!q}}

\newcommand{\qBinom}[2]{\genfrac{[}{]}{0pt}{1}{#1}{#2}}

\begin{document}

\rightline{LMU-TPW 2002-07}
\rightline{MPI-ThP/2002-43}

\vspace{3em}
\begin{center}
  
  {\Large{\bf Covariant Realization of Quantum Spaces as
      Star Products by Drinfeld Twists}}

\vspace{3em}

\textbf{Christian Blohmann}

\vspace{1em}
 
Ludwig-Maximilians-Universit\"at M\"unchen, Sektion Physik\\
Lehrstuhl Prof.\ Wess, Theresienstr.\ 37, D-80333 M\"unchen\\[1em]

Max-Planck-Institut f\"ur Physik, 
        F\"ohringer Ring 6, D-80805 M\"unchen\\[1em]

\end{center}

\vspace{1em}

\begin{abstract}
  Covariance of a quantum space with respect to a quantum enveloping
  algebra ties the deformation of the multiplication of the space
  algebra to the deformation of the coproduct of the enveloping
  algebra.  Since the deformation of the coproduct is governed by a
  Drinfeld twist, the same twist naturally defines a covariant star
  product on the commutative space.  However, this product is in
  general not associative and does not yield the quantum space. It is
  shown that there are certain Drinfeld twists which realize the
  associative product of the quantum plane, quantum Euclidean 4-space,
  and quantum Minkowski space. These twists are unique up to a central
  2-coboundary. The appropriate formal deformation of real structures
  of the quantum spaces is also expressed by these twists.
\end{abstract}


\section{Introduction}

Noncommutative geometries are described by replacing the commutative
algebra of functions on an ordinary space with a noncommutative
algebra. The noncommutativity is controlled by a perturbation
parameter, if it is a small deviation from ordinary geometry. The
algebraic aspects of such a perturbation can be detached from
questions of convergence and continuity by considering formal power
series. A noncommutative geometry is then described as a formal
deformation of a commutative algebra \cite{Gerstenhaber:1964} or by a
star product \cite{Bayen:1978a}. Such a description has attracted a
lot of attention lately, due to its application to string theory
\cite{Seiberg:1999} and to the construction of gauge theories on
noncommutative spaces \cite{Madore:2000b}.

The description of physical spacetime by an algebra alone would not
distinguish Minkowski space from, say, Euclidean 4-space, which
differs by the symmetry that acts on it. Deforming a space algebra
which transforms covariantly under a symmetry Lie group, will in
general break the symmetry. But there are deformations, where the
symmetry structure can be deformed together with the space, so that
covariance is preserved. Quantum spaces
\cite{Manin:1988,Faddeev:1990,Carow-Watamura:1990} are such a class of
deformations, carrying a covariant representation of the
Drinfeld-Jimbo deformation \cite{Drinfeld:1985,Jimbo:1985} of the
enveloping symmetry algebra.

The deformation of an enveloping algebra into a Drinfeld-Jimbo algebra
is well understood: As algebra over the ring of formal power series
the deformed algebra is isomorphic to the undeformed one. In fact, if
the Lie algebra is semisimple, it can be shown by cohomological
arguments that the enveloping algebra cannot be deformed at all
\cite{Gerstenhaber:1964}. It is only the Hopf structure which is truly
deformed. The deformed and the undeformed coproduct are non isomorphic
but related by inner automorphisms, called Drinfeld twists
\cite{Drinfeld:1989b}. Preserving covariance ties the deformation of
the enveloping algebra closely to the deformation of the space
algebra. Therefore, one ought to be able to use the knowledge about
the deformation of the symmetry in order to deform the space algebra
accordingly. Such an approach was suggested in \cite{Grosse:2001},
where a Drinfeld twist was used to realize a quantum space as star
product on the undeformed space algebra (for a similar approach see
\cite{Bonechi:1995}). By construction, such a star product is
covariant with respect to the action of the Drinfeld-Jimbo algebra.
However, the star product is in general not associative. The main
question of this article is: Is there a Drinfeld twist which
implements the associative product of a given quantum space? We will
give a positive answer for three particularly interesting cases: the
quantum plane, quantum Euclidean 4-space, and quantum Minkowski space.

In Sec.~\ref{sec:Twists} we review quantum spaces as covariant
deformations and relate them to star products defined by Drinfeld
twists \cite{Grosse:2001}. We recall the definition and properties of
Drinfeld twists in the framework of formal algebraic deformations.
Real structures of quantum spaces and of the according quantum
algebras are taken into account: Gerstenhaber's rigidity theorem for
algebras \cite{Gerstenhaber:1964} is extended to $*$-algebras
(Prop.~\ref{th:StarRigid}) and the deformation of space algebras is
extended to real structures (Prop.~\ref{th:StarModuleDeform}). In
Sec.~\ref{sec:StarProducts} we propose a general approach which
reduces the algebraic problem of finding a twist which implements the
multiplication of a quantum space to a representation theoretic
problem. This works well for cases where the representation theory of
the symmetry quantum algebra is well understood: We determine the
basis which reduces the quantum plane, quantum Euclidean 4-space, and
quantum Minkowski space as module into its irreducible highest weight
subrepresentations and calculate the multiplication map with respect
to this basis. Comparing the multiplication maps with the
representations of the twists leads to the main result: There are
Drinfeld twists which realize the quantum plane
(Prop.~\ref{th:PlaneF}), quantum Euclidean 4-space
(Prop.~\ref{th:EuclidF}), and quantum Minkowski space
(Prop.~\ref{th:MinkowskiF}) as covariant star products. These twists
are unique up to a central 2-coboundary.

Throughout this article we assume that $\mathfrak{g}$ is a semisimple
Lie algebra, denoting its enveloping algebra by $\Ug$. An element $u
\in \Ug^{\otimes (n+1)}$ is called $\mathfrak{g}$-invariant if it
commutes with the $n$-fold coproduct $\Delta^{(n)}(g) := g \otimes
1^{\otimes n} + 1 \otimes g \otimes 1^{\otimes (n-1)} + \ldots +
1^{\otimes n} \otimes g$ of all $g \in \mathfrak{g}$. The formal
perturbation parameter is $\h$, the completion of a complex vector
space or algebra $A$ with respect to the $\h$-adic topology is
$A[[\h]]$. The topological tensor product $\hat{\otimes}$ of two free
$\h$-adic vector spaces or algebras is avoided by identifying $A[[\h]]
\hat{\otimes} A'[[\h]] \equiv (A \otimes A')[[\h]]$. The $\h$-adic
Drinfeld-Jimbo deformation \cite{Drinfeld:1985,Jimbo:1985} of $\Ug$ is
denoted by $\Uhg$. The equality of two elements $a,a' \in A[[\h]]$
modulo $\h^n$ will be written in Landau notation as $a = a' +
\Ocal(\h^n)$. Recall that if $a = 1 + \Ocal(\h)$ then $a$ is
invertible and its square root with $\sqrt{a} = 1 +\Ocal(\h)$ is
defined and unique in $A[[\h]]$ (see e.g. \cite{Kassel}). The
symmetric $\h$-adic quantum number is defined by $[n]:= (\E^{\h
  n}-\E^{-\h n})(\E^{\h}-\E^{-\h})^{-1}$ and for natural $n$ the
quantum factorial by $[n]! := [1]\cdot[2]\cdots[n]$.

\section{Quantum Spaces and Drinfeld Twists}
\label{sec:Twists}

\subsection{Quantum Spaces}

\label{sec:QuantumSpaces}

Let $\mathfrak{g}$ be the Lie algebra of the symmetry group of a space
and $\Xcal$ be the function algebra of this space. The elements $g \in
\mathfrak{g}$ of the Lie algebra act on $\Xcal$ as derivations, $g \tr
xy = (g\tr x)y + x(g \tr y)$ for $x,y \in \Xcal$. A generalized way of
writing this is
\begin{equation}
\label{eq:ModuleAlgebra}
  g \tr xy = (g_{(1)} \tr x)(g_{(2)} \tr y)
\end{equation}
for all $g \in \Ug$, where we introduce the coproduct of an enveloping
algebra by $g_{(1)} \otimes g_{(2)} \equiv \Delta(g) := g \otimes 1 +
1 \otimes g$ on the generators $g \in \mathfrak{g}$ and extend it to a
homomorphism $\Delta: \Ug \rightarrow \Ug \otimes \Ug$ on the
enveloping algebra.

Using the multiplication map $\mu : \Xcal \otimes \Xcal \rightarrow
\Xcal$, $\mu(x \otimes y) := xy$ of $\Xcal$ we can
write~\eqref{eq:ModuleAlgebra} as
\begin{equation}
\label{eq:ModuleAlgebra2}
  g \tr \mu(x \otimes y) = \mu(\Delta(g) \tr [x \otimes y])\,,
\end{equation}
the condition for the product $\mu$ to be covariant with respect to
the action of $\Ug$, which is meaningful not only for $\Ug$ but for
any Hopf algebra. In mathematical terminology, an algebra $\Xcal$
which carries a representation of some Hopf algebra $H$ such that
Eq.~\eqref{eq:ModuleAlgebra} holds is called an $H$-module algebra. We
will also call it an $H$-covariant space. The covariant spaces or
module algebras of the quantum enveloping algebras $\Uhg$ are called
quantum spaces.

As algebra, $\Uhg$ for semisimple $\mathfrak{g}$ is not a true
deformation of $\Ug$, because $\Uhg$ and $\Ug[[\h]]$ are isomorphic as
algebras (see Sec.~\ref{sec:Drinfeld1}). This means that every
$\h$-adic space algebra $\Xcal$ which is a $\Ug[[\h]]$-module is also
a $\Uhg$-module and vice versa. In other words, if we considered only
the algebra structure of $\Ug$, then there would be no need to replace
a commutative space with a noncommutative one when passing from the
symmetry algebra to its quantum deformation. It is the Hopf structure
which is deformed in an essential way. That is, identifying $\Uhg$ and
$\Ug[[\h]]$ as isomorphic algebras, we can view the quantum
deformation as deformation of the Hopf structure of $\Ug[[\h]]$,
$\Delta \rightarrow \Delta_\h$, $\varepsilon \rightarrow
\varepsilon_\h$, $S \rightarrow S_\h$. Since the covariance
condition~\eqref{eq:ModuleAlgebra} ties the algebra structure of a
space $\Xcal$ to the coproduct, the multiplication map of the space
must be deformed, $\mu \rightarrow \mu_\h$, along with the deformation
of the coproduct, $\Delta \rightarrow \Delta_\h$, if the covariance is
to be preserved. Conversely, deforming the multiplication map of a
covariant space, the coproduct of the symmetry algebra must be
deformed accordingly,
\begin{equation}
\label{eq:BothDeform}
\begin{CD}
  g \tr xy = (g_{(1)} \tr x)(g_{(2)} \tr y) 
    @>{\mu \rightarrow \mu_\h}>{\Delta \rightarrow \Delta_\h}>
  g \tr (x \star y) = (g_{(1_\h)} \tr x)\star (g_{(2_\h)} \tr y) \,,
\end{CD}
\end{equation}
where $\Delta_\h(g) = g_{(1_\h)} \otimes g_{(2_\h)}$ and $\mu_\h(x
\otimes y) = x \star y$. While there may be a large class of
deformations which are covariant in this sense, we will restrict our
attention to quantum spaces.

\subsection{Star Products by Drinfeld Twists}
\label{sec:Quest}

In the case of quantum spaces, the deformed coproduct belongs to the
Drinfeld-Jimbo deformation $\Uhg \cong (\Ug[[\h]], \Delta_\h,
\varepsilon_\h, S_\h)$. Drinfeld has observed
(Theorem~\ref{th:Drinfeld2}) that $\Delta_\h$ is related to the
undeformed coproduct $\Delta$ by an inner automorphism. That is, there
is an invertible element $\F \in (\Ug \otimes \Ug)[[\h]]$ with $\F = 1
+ \Ocal(\h)$, called Drinfeld twist, such that
\begin{equation}
\label{eq:StarProd4}
  \Delta_\h(g) := \F\Delta(g)\F^{-1} \,.
\end{equation}
Comparing the covariance condition of the deformed multiplication,
\begin{equation}
\label{eq:StarProd3}
  g \tr \mu_\h(x \otimes y) = \mu_\h(\Delta_\h(g)\tr [x \otimes y])
  = \mu_\h(\F\Delta(g)\F^{-1} \tr [x \otimes y])
\end{equation}
with the covariance property~\eqref{eq:ModuleAlgebra2} of the
undeformed product, we see that Eq.~\eqref{eq:StarProd3} is naturally
satisfied if we define the deformed product by \cite{Grosse:2001}
\begin{equation}
\label{eq:StarProd2}
  \mu_\h(x \otimes y) := \mu(\F^{-1} \tr [x \otimes y])
  \quad \Leftrightarrow \quad
  x \star y := (\F^{-1}_{[1]} \tr x)(\F^{-1}_{[2]} \tr y)\,,
\end{equation}
where we suppress in a Sweedler like notation the summation of $\sum_i
\F_{1i} \otimes \F_{2i} \equiv \F_{[1]} \otimes \F_{[2]}$. Since the
elements of the Lie algebra $\mathfrak{g}$ act on the undeformed space
algebra $\Xcal$ as derivations, $\F^{-1}$ acts as $\h$-adic
differential operator on $\Xcal \otimes \Xcal$. Hence, writing out the
$\h$-adic sum of $\F^{-1} = 1 + \sum_k \h^k \F_k^{-1}$ we can define
the bidifferential operators
\begin{equation}
\label{eq:Bdef}
  B_k(x \otimes y) := \mu(\F_k^{-1} \tr [x \otimes y])
  = (\F^{-1}_{k[1]} \tr x)(\F^{-1}_{k[2]} \tr y)\,,
\end{equation}
such that the star product~\eqref{eq:StarProd2} can be written in the
more familiar form
\begin{equation}
\label{eq:StarProd1}
  x \star y := xy + \h B_1(x,y) + \h^2 B_2(x,y) + \ldots 
\end{equation}
Even though the twist $\F$ yields by Eq.~\eqref{eq:StarProd4} a
coassociative coproduct, Eq.~\eqref{eq:StarProd2} will in general not
define an associative product. The associativity condition for
$\mu_\h$ reads
\begin{equation}
\label{eq:ProdAssoc}
\begin{split}
  (x \star y) \star z
  &= (\F^{-1}_{[1](1)}\F^{-1}_{[1']} \tr x)
    (\F^{-1}_{[1](2)}\F^{-1}_{[2']} \tr y)
    (\F^{-1}_{[2]} \tr z ) \\
  &= (\F^{-1}_{[1]} \tr x )
     (\F^{-1}_{[2](1)}\F^{-1}_{[1']} \tr y)
     (\F^{-1}_{[2](2)}\F^{-1}_{[2']} \tr z)
  = x \star (y \star z) \,,
\end{split}  
\end{equation}
for all $x,y,z \in \Xcal$. Defining the Drinfeld coassociator
\begin{equation}
\label{eq:CoassDef}
\Phi := (\Delta \otimes \id)(\F^{-1})\, (\F^{-1} \otimes 1)
  \,(1 \otimes \F) \,(\id \otimes \Delta)(\F) \,,
\end{equation}
the associativity condition~\eqref{eq:ProdAssoc} can be written as
\begin{equation}
\label{eq:Quest1}
  (\Phi_{[1]} \tr x) (\Phi_{[2]} \tr y)(\Phi_{[3]} \tr z)
  = xyz \,.
\end{equation}
The obvious question to ask is: For a given $\Uhg$-covariant quantum
space, is there a Drinfeld twist $\F$ which yields by
Eq.~\eqref{eq:StarProd2} the associative product of the quantum space?
We do not attempt to answer this question in its generality. Instead,
we will consider some prototypical and physically important cases: the
quantum plane, quantum Euclidean 4-space, and quantum Minkowski space.

\subsection{Drinfeld Twists of Quantum Enveloping Algebras}
\label{sec:Drinfeld1}

For the reader's convenience we gather in this section some well known
results on formal deformations of algebras and Hopf algebras,
essentially due to Gerstenhaber \cite{Gerstenhaber:1964} and Drinfeld
\cite{Drinfeld:1989,Drinfeld:1989b}.

An $\h$-adic algebra $A'$ is called a deformation of an algebra $A$ if
$A'/\h A'$ and $A$ are isomorphic as algebras. Analogously, an
$\h$-adic Hopf algebra $H'$ is called a deformation of a Hopf algebra
$H$ if $H'/ \h H'$ and $H$ are isomorphic as Hopf algebras. Recall
that $\Ug$ is a Hopf algebra with the canonical Lie Hopf structure
defined on the generators $g\in \mathfrak{g}$ as coproduct $\Delta(g)
= g \otimes 1 + 1 \otimes g$, counit $\varepsilon(g) = 0$, and
antipode $S(g) = -g$. The Drinfeld-Jimbo algebra $\Uhg$ is a
deformation of this Hopf algebra $\Ug$. This can be seen by developing
the commutation relations and the Hopf structure of $\Uhg$ as formal
power series in $\h$ and keeping only the zeroth order terms, which
yields the commutation relations and the Lie Hopf structure of $\Ug$.

Gerstenhaber has shown \cite{Gerstenhaber:1964} that whenever the
second Hochschild cohomology of $A$ with coefficients in $A$ is zero,
$H^2(A,A) = 0$, then all deformations of $A$ are trivial up to
isomorphism. That is, any deformation $A'$ of $A$ is isomorphic to the
$\h$-adic completion of the undeformed algebra, $A' \cong A[[\h]]$.
Algebras with this property are called rigid. The second Whitehead
lemma states that the second Lie algebra cohomology of a semisimple
Lie algebra $\mathfrak{g}$ and, hence, the second Hochschild
cohomology of its enveloping algebra is zero. Therefore, the
enveloping algebra $\Ug$ of a semisimple Lie algebra $\mathfrak{g}$ is
rigid. In particular, there is an isomorphism of algebras $\alpha:
\Uhg \rightarrow \Ug[[\h]]$, by which the the Hopf structure
$\Delta'$, $\varepsilon'$, $S'$ of $\Uhg$ can be transfered to
$\Ug[[\h]]$,
\begin{equation}
\label{eq:HopfDeform}
  \Delta_\h := (\alpha\otimes\alpha)
    \circ \Delta' \circ\alpha^{-1} \,,\quad
  \varepsilon_\h := \varepsilon' \circ \alpha^{-1} \,,\quad
  S_\h := \alpha \circ S' \circ \alpha^{-1} \,,
\end{equation}
such that $\alpha$ becomes an isomorphism of Hopf algebras from $\Uhg$
to $\Ug[[\h]]$ with this deformed Hopf structure. Let $\alpha'$ be
another such isomorphism and $\Delta'_\h$, $\varepsilon'_\h$, $S'_\h$
be defined as in Eq.~\eqref{eq:HopfDeform} with $\alpha'$ instead of
$\alpha$. Then $\alpha'$ is an isomorphism of Hopf algebras from
$\Uhg$ to $\Ug[[\h]]$ with the primed Hopf structure,
\begin{equation}
  (\Ug[[\h]], \Delta_\h, \varepsilon_\h, S_\h)
  \stackrel{\alpha}{\longleftarrow} \Uhg
  \stackrel{\alpha'}{\longrightarrow}
  (\Ug[[\h]], \Delta'_\h, \varepsilon'_\h, S'_\h) \,,
\end{equation}
hence, $\alpha'\circ \alpha^{-1}$ is an isomorphism of Hopf algebras.
We conclude that, while the Hopf structure Eq.~\eqref{eq:HopfDeform}
may depend on the isomorphism $\alpha$, it is unique up to an
isomorphism of Hopf algebras. 

As a consequence of the first Whitehead lemma, the first Hochschild
cohomology of the enveloping algebra $\Ug$ of a semisimple Lie algebra
is zero. This implies, that the two homomorphisms $\Delta$ and
$\Delta_\h$ from $\Ug[[\h]]$ to $(\Ug \otimes \Ug)[[\h]]$ with
$\Delta_\h = \Delta + \Ocal(\h)$ are related by an inner automorphism,
as it was observed by Drinfeld \cite{Drinfeld:1989,Drinfeld:1989b}.

\begin{Theorem}
\label{th:Drinfeld2}
Let $\mathfrak{g}$ be a semisimple Lie algebra, and let $\Delta_\h$ be
defined as in Eq.~\eqref{eq:HopfDeform}. Then there is an invertible
element $\F \in \bigl(\Ug \otimes \Ug\bigr)[[\h]]$ such that
$\Delta_\h(g) = \F \Delta(g) \F^{-1}$, which is called a Drinfeld
twist from $\Delta$ to $\Delta_\h$.
\end{Theorem}

On the first sight Theorem~\ref{th:Drinfeld2} only relates the
coproducts. It turns out that the twist of the coproduct relates
counit and antipode, as well.

\begin{Corollary}
\label{th:HopfTwist}
Let $\F$ be a Drinfeld twist from $\Delta$ to $\Delta_\h$ as in
Theorem~\ref{th:Drinfeld2}.
\begin{itemize}
  
\item[\textup{(i)}] If $\F'$ is another Drinfeld twist, then $\F^{-1}
  \F'$ is invertible and $\mathfrak{g}$-invariant.  Conversely, let
  $\Tcal \in \bigl(\Ug \otimes \Ug\bigr)[[\h]]$ be invertible and
  $\mathfrak{g}$-invariant.  Then $\F \Tcal$ is a Drinfeld twist.

\item[\textup{(ii)}] $\varepsilon_\h = \varepsilon$
  
\item[\textup{(iii)}] There is a twist $\F$ such that
  $\varepsilon(\F_{[1]})\F_{[2]} = 1 = \F_{[1]}\varepsilon(\F_{[2]})$,
  which implies $\F = 1 + \Ocal(\h)$.  Twists with this property are
  called counital.
  
\item[\textup{(iv)}] The two elements of $\Ug[[\h]]$ defined as
  \begin{equation}
  \label{eq:sigmaDefine}
    \sigma_1^{-1} := S(\F_{[1]}^{-1})\F_{[2]}^{-1} \,,\qquad
    \sigma_2 := \F_{[1]} S(\F_{[2]}) 
  \end{equation}
  are invertible, $\sigma_1^{-1}\sigma_2 = \sigma_2\sigma_1^{-1}$ is
  central in $\Ug[[\h]]$, and $\sigma_1 S(g) \sigma_1^{-1} = S_\h(g) =
  \sigma_2 S(g) \sigma_2^{-1}$ for all $g\in\Ug[[\h]]$.
  
\item[\textup{(v)}] The coassociator $\Phi$ defined as in
  Eq.~\eqref{eq:CoassDef} is $\mathfrak{g}$-invariant.
  
\item[\textup{(vi)}] The deformed Hopf structure on $\Ug[[\h]]$ is
  isomorphic to the undeformed one if and only if there is a Drinfeld
  twist of the form $\F = (u \otimes u)\Delta u^{-1}$ for some
  invertible $u \in \Ug[[\h]]$.

\end{itemize}
\end{Corollary}

A proof can be found in the appendix. The multiplication of a twist
with a $\mathfrak{g}$-invariant element as in (i) is sometimes called
a gauge transformation of the twist. Counitality (iii) is often part
of the definition of twists. Therefore, we will assume from now on
that all Drinfeld twists are counital. It can be shown that there are
twists for which the elements $\sigma_1$ and $\sigma_2$ of (iv) are
equal \cite{Fiore:1997}. However, this is not the case for all twists.
For example, assume that $\sigma_1 = \sigma_2$ for some twist $\F$,
and assume that there is a cubic Casimir $c \in \Ug$ with $S(c) = -c$.
Then $\F' := \F[1 \otimes (1 + \h c)]$ is another Drinfeld twist for
which $\sigma^{\prime -1}_1\sigma'_2 = \sigma_1^{-1}\sigma_2 (1+\h
c)^{-1}S(1+\h c) = (1+\h c)^{-1}(1 -\h c) \neq 1$. Part (vi) of the
corollary applies also to the case where Eq.~\eqref{eq:HopfDeform}
defines for two different isomorphism $\alpha$ and $\alpha'$ two
different coproducts $\Delta_\h$ and $\Delta'_\h$. Here, the
automorphism which relates the two coproducts is $\beta := \alpha'
\alpha^{-1}$. Note that elements of the form $(u \otimes u)\Delta
u^{-1}$ are 2-coboundaries in the sense of \cite{Majid}. Since
$\Delta$ is cocommutative, 2-coboundaries are symmetric. Hence, the
twist~\eqref{eq:StarProd4} of $\Delta$ by a coboundary yields a
cocommutative coproduct. The coproduct of the Drinfeld-Jimbo
deformation is not cocommutative, so it cannot be isomorphic to the
undeformed, cocommutative coproduct.

\subsection{Real Forms of Enveloping Algebras}

Lie groups are usually viewed as real manifolds, even though they may
be naturally defined as complex matrix groups. For example, the
universal covering of the Lorentz group $SL(2,\mathbb{C})$ is viewed
as real 6-parameter Lie group, the generators being the three
rotations and the three boosts. When considering the complexification
$\tilde{\mathfrak{g}} := \mathbb{C} \otimes_{\mathbb{R}} \mathfrak{g}$
of a real Lie algebra, we have to keep in mind that non-isomorphic
real Lie algebras can have the same complexification. For example
$\mathrm{su}_2$ and $\mathrm{sl}_2(\mathbb{R})$ have the same
complexification $A_1$. A practical method to remember the real Lie
algebra which a complexification comes from is to observe that for any
real Lie algebra $\mathfrak{g}$ there is an antihomomorphism $*$
defined as $g^* := -g$ for all $g \in \mathfrak{g}$, which can be
extended to a conjugate linear antihomomorphism on the
complexification by $(\alpha \otimes_{\mathbb{R}} g)^* :=
\overline{\alpha} \otimes_{\mathbb{R}} g^*$ for all $\alpha \in
\mathbb{C}$, $g \in \mathfrak{g}$. This defines a $*$-structure on
$\tilde{\mathfrak{g}}$, that is, a conjugate linear antihomomorphism,
which is an involution, $*^2 = \id$. Given the $*$-structure we can
reconstruct the real Lie algebra as being generated by all elements of
the form $(g - g^*) \in \tilde{\mathfrak{g}}$.

The identification of real forms of a complex Lie algebra with
$*$-structures can be extended to the enveloping Hopf algebra $\Ug$.
In general, a $*$-structure on a Hopf algebra is a conjugate linear
antihomomorphism $*$, which is an involution and a bialgebra
homomorphism, $\Delta(g^*) = (\Delta g)^{* \otimes *}$. Algebras and
Hopf algebras with such a $*$-structure are called $*$-algebras and
Hopf $*$-algebras, respectively. Finally, if $H$ is a Hopf $*$-algebra
and $\Xcal$ an $H$-module algebra with a $*$-structure, such that in
addition to Eq.~\eqref{eq:ModuleAlgebra} we have
\begin{equation}
\label{eq:StarModule}
  (g \tr x)^* = (Sg)^* \tr x^*
\end{equation}
for all $g \in H$, $x \in \Xcal$, then $\Xcal$ is called an $H$-module
$*$-algebra, or a real $H$-covariant space.

Analogously as for algebras, an $\h$-adic $*$-algebra $A'$ is called a
deformation of a $*$-algebra $A$ if $A'/ \h A'$ and $A$ are isomorphic
as $*$-algebras. $A$ is called rigid as $*$-algebra if for any
deformation $A'$ of $A$ the $\h$-adic completion $A[[\h]]$ and $A'$
are isomorphic as $*$-algebras.

\begin{Proposition}
\label{th:StarRigid}
Let $A$ be a $*$-algebra with zero first and second Hochschild
cohomology, $H^1(A,A) = H^2(A,A) = 0$. Then $A$ is rigid as
$*$-algebra.
\end{Proposition}

\begin{proof}
  Let $(A',*')$ be a deformation of $(A,*)$ as $*$-algebra. $H^2(A,A)
  = 0$ implies that $A$ is rigid as algebra, so there is an
  isomorphism of algebras $\beta : A' \rightarrow A[[\h]]$. Define
  $*_\beta := \beta \circ *' \circ \beta^{-1}$.  Because $*'$ is a
  deformation of $*$, we have $*_\beta = * + \Ocal(\h)$, and thus, as
  $*^2 = \id$, $*_\beta \circ * = \id + \Ocal(\h)$. $H^1(A,A) = 0$
  implies that this algebra automorphisms is inner, $(g)^{*_\beta
    \circ *} = u g u^{-1}$ for some invertible $u \in A[[\h]]$. Thus,
  $g^{*_\beta} = (u^*)^{-1} g^* u^*$ for all $g \in A[[\h]]$. Since
  $*_\beta^2 = \id$, we have $(g^{*_\beta})^{*_\beta} = g =
  (u^*)^{-1}u g u^{-1} u^*$, so $u^{-1} u^* = u^* u^{-1}$ is central.
  Since $*_\beta$ is a deformation of $*$ we can choose $u$ such that
  $u = 1 + \Ocal(\h)$ and, thus, $u^{-1} u^* = 1 + \Ocal(\h)$, so the
  square roots of $u^*$ and $u^{-1} u^*$ are defined and invertible.
  Define an isomorphism of algebras $\alpha: A' \rightarrow A[[\h]]$
  by
  \begin{equation}
    \alpha(g) := (u^*)^{\frac{1}{2}} \, \beta(g)\, (u^*)^{-\frac{1}{2}} \,.
  \end{equation}
  Using that the square root of a central element is central, we get
  \begin{equation}
  \begin{split}
    g^{*_\alpha} &:= (\alpha \circ *' \circ \alpha^{-1})(g) =
    (u^*)^{\frac{1}{2}} \, (\beta \circ *' \circ \beta^{-1}) \bigl[
    (u^*)^{-\frac{1}{2}} \, g\, (u^*)^{\frac{1}{2}}
    \bigr] \, (u^*)^{-\frac{1}{2}} \\
    &= (u^*)^{\frac{1}{2}} (u^*)^{-1} (u)^{\frac{1}{2}} g^*
    (u)^{-\frac{1}{2}} u^* (u^*)^{-\frac{1}{2}} = (u^{-1}
    u^*)^{-\frac{1}{2}} \, g^* \, (u^{-1} u^*)^{\frac{1}{2}} = g^* \,,
  \end{split}
  \end{equation}
  which shows that $\alpha:(A',*') \rightarrow (A[[\h]], *)$ is an
  isomorphism of $*$-algebras.
\end{proof}

This proposition applies in particular to the real forms of enveloping
algebras of semisimple Lie algebras. In fact, if $*$ is a
$*$-structure on $\Ug$ and $*'$ is a $*$-structure on $\Uhg$, such
that $(\Uhg,*')$ is a deformation of $(\Ug, *)$ as Hopf $*$-algebra,
we can (and shall) always use an isomorphism of $*$-algebras
$\alpha:(\Uhg,*') \rightarrow(\Ug, *)$ to transfer the Hopf structure
by Eqs.~\eqref{eq:HopfDeform}. In this case there are twists with
particularly interesting properties with respect to the $*$-structure:
We will call a twist unitary and real, respectively, if
\begin{subequations}
\label{eq:Forthogonal}
\begin{alignat}{2}
\label{eq:Funitary}
  &\text{unitary:} \quad &(* \otimes *)(\F) &= \F^{-1} \\
\label{eq:Freal}
  &\text{real:}\quad &(* \otimes *)(\F) &= (S \otimes S)(\F_{21}) \,.
\end{alignat}
\end{subequations}
A twist which is both, unitary and real, is called orthogonal.
\begin{Proposition}
\label{th:Forth}
There is an orthogonal Drinfeld twist from $\Delta$ to $\Delta_\h$.
\end{Proposition}
\begin{proof}
  In Theorem~4.1 of \cite{Donin:1997} it was shown that there is
  always a twist with $(S \otimes S)(\F) = \F_{21}^{-1}$. By
  assumption, the $*$-structure is a homomorphism of coalgebras for
  both, $\Delta$ and $\Delta_\h$. Thus, $(* \otimes * )(\F^{-1})
  \equiv (\F^{-1})^*$ is also a twist, so $\Tcal := \F^{-1}(\F^{-1})^*
  = 1 + \Ocal(\h)$ is $\mathfrak{g}$-invariant and $\F' := \F
  \sqrt{\Tcal}$ is another twist. It is easy to check that $\F'$ is
  unitary and real.
\end{proof}

\subsection{Real Structures on Star Products}

The quantum spaces we want to consider here, the quantum plane and
quantum Minkowski space, possess a real structure which is covariant
with respect to a real form of the symmetry algebra. That is,
$(\Ug,*)$ is a Hopf $*$-algebra and $(\Xcal, \mu, *)$ is a module
$*$-algebra, the action of $\Ug$ on $\Xcal$ and the $*$-structures
satisfying Eq.~\eqref{eq:StarModule}. For a general twist
Eq.~\eqref{eq:StarProd2} will not define a multiplication which is
compatible with any real structure. Assume that we deform the product
$\mu \rightarrow \mu_\h$ of the real covariant space $(\Xcal, \mu, *)$
but not the $*$-structure. Then
\begin{equation}
\label{eq:antihom}
\begin{split}
   (x \star y)^*
   &= [(\F_{[1]}^{-1} \tr x)(\F_{[2]}^{-1} \tr y)]^*
   = (\F_{[2]}^{-1} \tr y)^*(\F_{[1]}^{-1} \tr x)^* \\
   &= ([S\F_{[2]}^{-1}]^* \tr y^*) ([S\F_{[1]}^{-1}]^* \tr x^*) \\
   &= (\F_{[2]}[S\F_{[2]}^{-1}]^* \tr y)^* \star
   (\F_{[1]}[S\F_{[1]}^{-1}]^* \tr x)^* \,,
\end{split}
\end{equation}
which shows that the undeformed $*$-structure is an antihomomorphism
with respect to the deformed product, if $\F$ is chosen to be real in
the sense of Eq.~\eqref{eq:Freal}, which was proved to be possible in
Prop.~\ref{th:Forth}. However, the undeformed $*$-structure of $\Xcal$
will in general not satisfy the module $*$-algebra
property~\eqref{eq:StarModule}. Thus, while Prop.~\ref{th:StarRigid}
shows that the undeformed and deformed $*$-structures of the symmetry
algebra $\Ug[[\h]]$ can (and shall) be chosen to coincide, this is not
possible for the $*$-structures of the deformed and undeformed space
algebras. However, there is a unique element $\sigma$ of the symmetry
algebra which mediates the deformation $* \rightarrow *_\h$ of the
real structure of $\Xcal$ by $x^{*_\h} := \sigma^{-1} \tr x^*$. This
element is characterized by
\begin{Proposition}
\label{th:sigmaunique}
There is a unique element $\sigma \in \Ug[[\h]]$ such that
  \begin{equation}
  \label{eq:sigmaDef}
    \sigma = 1 + \Ocal(\h) \,,\qquad
    S_\h(g) = \sigma(Sg)\sigma^{-1}\,, \qquad
    \Delta_\h(\sigma) = \sigma \otimes \sigma \,.
  \end{equation}
Moreover, $\sigma^* = \sigma$.
\end{Proposition}
\begin{proof}
  Let $\alpha: \Uhg \rightarrow \Ug[[\h]]$ be the isomorphism of
  algebras which is used to define the deformed Hopf structure by
  Eqs.~\eqref{eq:HopfDeform}. Let $\R$ be a universal $\R$-matrix of
  $\Uhg$ and
  \begin{equation}
  \label{eq:RDeform}
    \R_\h := (\alpha \otimes \alpha)(\R) \,,
  \end{equation}
  such that $\R_\h$ becomes a universal $\R$-matrix with respect to
  $\Delta_\h$.  It was shown in Prop.~3.16 of \cite{Drinfeld:1989b}
  and Theorem~4.1 of \cite{Donin:1997} that there is a Drinfeld twist
  $\F$ from $\Delta$ to $\Delta_\h$ such that
  \begin{equation}
  \label{eq:RFRel}
    \R_\h = \F_{21} \E^{\frac{\h}{2}
    (\Delta(C) - C \otimes 1 - 1 \otimes C) } \F^{-1}
    \qquad\text{and}\qquad (S \otimes S)(\F) = \F_{21}^{-1} \,,  
  \end{equation}
  with the canonical quadratic Casimir $C := g_i g_j K^{ij}$, where
  $\{g_i\}$ is a basis of $\mathfrak{g}$, $K_{ij} := \mathrm{tr}(
  \mathrm{ad}\,g_i\, \mathrm{ad}\,g_j)$ is the Killing metric, and
  $K^{ij}$ its inverse, $K^{ij}K_{jk} = \delta^i_k$.
  
  Let $\sigma_1^{-1}$ and $\sigma_2$ be defined for such a twist as in
  Eq.~\eqref{eq:sigmaDefine}, and let $\sigma := \sqrt{\sigma_1
    \sigma_2}$.  Since $\sigma_1^{-1} \sigma_2$ is central, so is its
  square root.  Hence, $S_\h(g) = \sigma_1(Sg) \sigma_1^{-1} =
  \sigma_1 \sqrt{\sigma_1^{-1} \sigma_2} (Sg) \sqrt{\sigma_2^{-1}
    \sigma_1} \sigma_1^{-1} = \sigma (Sg) \sigma^{-1}$.
  
  From the second equation of~\eqref{eq:RFRel} we deduce $S(\sigma_2)
  = \sigma_1^{-1}$, so $S(\sigma) = \sigma^{-1}$.  From the first
  equation of~\eqref{eq:RFRel} we compute $u:=
  S_\h(\R_{\h[2]})\R_{\h[2]} = \sigma_1 S(\sigma_1^{-1}) q^{-C} =
  \sigma_1 \sigma_2 q^{-C} = \sigma^2 q^{-C}$.  From the properties of
  universal $\R$-matrices it follows that $uS_\h(u^{-1}) = \sigma^4$
  is group-like with respect to $\Delta_\h$, $\Delta_\h(\sigma_1^4) =
  \sigma_1^4 \otimes \sigma_1^4$. Since $\sigma^4$ is group-like, its
  fourth root $\sigma$ is group-like, as well.
  
  Now let $\sigma' = 1 + \Ocal(\h)$ be another element with $S_\h(g) =
  \sigma'(Sg) \sigma^{\prime -1}$ which is group-like with respect to
  $\Delta_\h$. Then $z := \sigma'\sigma^{-1}$ is central and
  group-like, $\Delta_\h(z) = z \otimes z = \F^{-1}(z \otimes z) \F =
  \Delta(z)$, with respect to both coproducts.  Since $z = 1 +
  \Ocal(\h)$ there is an $a$ such that $z = \E^{\h a}$. Since $z$ is
  group-like and central $a$ must be primitive, $\Delta(a) = a \otimes
  1 + 1 \otimes a$, and central in every order of $a = \sum_i \h^i
  a_i$.  Being primitive, the $a_i$ are elements of the Lie algebra,
  $a_i \in \mathfrak{g} \subset \Ug$. Since $\mathfrak{g}$ is
  semisimple it does not contain nontrivial central elements, so $a_i
  = 0$ for all $i$, that is, $a = 0$. Hence, $z = \sigma'\sigma^{-1} =
  1$.
  
  Clearly, $\sigma^* = 1 + \Ocal(\h)$.  From $S_\h = * \circ S_\h^{-1}
  \circ *$ it follows that $S_\h(g) = \sigma^* (Sg) (\sigma^*)^{-1}$.
  Since $\sigma$ is group-like with respect to $\Delta_\h$, so is
  $\sigma^*$. By uniqueness of an element with
  properties~\eqref{eq:sigmaDef} we conclude that $\sigma^* = \sigma$.
\end{proof}

Now we can show, that $\sigma$ realizes the deformation of the real
structure in the promised manner:
\begin{Proposition}
\label{th:StarModuleDeform}
Let $\F$ be a twist from $\Delta$ to $\Delta_\h$ which is real in the
sense of Eq.~\eqref{eq:Freal} and let $(\Xcal, \mu, *)$ be an
$\h$-adic module $*$-algebra of $(\Ug[[\h]], \Delta, \varepsilon, S,
*)$. Define $\mu_\h$ as in Eq.~\eqref{eq:StarProd2} and $*_\h : \Xcal
\rightarrow \Xcal$ by
  \begin{equation}
  \label{eq:StarDeform}
    x^{*_\h} := \sigma^{-1} \tr x^*  
  \end{equation}
  for all $x \in \Xcal$, where $\sigma$ is the unique element of
  Prop.~\ref{th:sigmaunique}. Then $(\Xcal, \mu_\h, *_\h)$ is a module
  $*$-algebra of $(\Ug[[\h]], \Delta_\h, \varepsilon_\h, S_\h, *)$
\end{Proposition}
\begin{proof}
  By construction, $*_\h$ is conjugate linear.  We have to show that
  $*_\h$ is an algebra antihomomorphism. Writing $\mu_\h(x \otimes y)
  = x \star y$, we get
  \begin{align}
     (x \star y)^{*_\h}
     &= \sigma^{-1} \tr [(\F_{[1]}^{-1} \tr x)(\F_{[2]}^{-1} \tr y)]^*
     = \sigma^{-1} \tr (\F_{[2]}^{-1} \tr y)^*(\F_{[1]}^{-1} \tr x)^* \notag\\
     &= \sigma^{-1} \tr ([S\F_{[2]}^{-1}]^* \tr y^*)
                   ([S\F_{[1]}^{-1}]^* \tr x^*)     
     = \sigma^{-1} \tr (\F_{[1]}^{-1} \tr y^*)(\F_{[2]}^{-1} \tr x^*) \notag\\
     &= (\sigma^{-1}_{(1)}\F_{[1]}^{-1} \tr y^*)
       (\sigma^{-1}_{(2)}\F_{[2]}^{-1} \tr x^*)
     = (\F_{[1]}^{-1}\sigma^{-1}_{(1_\h)} \tr y^*)
       (\F_{[2]}^{-1}\sigma^{-1}_{(2_\h)} \tr x^*) \notag\\
     &= (\F_{[1]}^{-1}\sigma^{-1} \tr y^*)
        (\F_{[2]}^{-1}\sigma^{-1} \tr x^*)
     = (\sigma^{-1} \tr y^*)\star (\sigma^{-1} \tr x^*) \notag \\
     &= (y^{*_\h}) \star (x^{*_\h}) \,,
  \end{align}
  where we have used the module $*$-algebra
  condition~\eqref{eq:StarModule}, the assumed
  reality~\eqref{eq:Freal} of $\F$, and that $\sigma$ is group-like
  with respect to $\Delta_\h$. Since
  \begin{equation}
    (x^{*_\h})^{*_\h} = \sigma^{-1} \tr (\sigma^{-1} \tr x^*)^*
    = \sigma^{-1}(S\sigma^{-1})^* \tr x = \sigma^{-1}\sigma^* \tr x = x
  \end{equation}
  $*_\h$ is an involution. Finally, the module $*$-algebra
  condition~\eqref{eq:StarModule} holds,
  \begin{equation}
  \begin{split}
    (g \tr x)^{*_\h} &= \sigma^{-1} \tr (g \tr x)^*
    = \sigma^{-1} (Sg)^* \tr x^*
    = (\sigma^{-1})^* (Sg)^* \tr x^* \\
    &=  (S_\h g)^* (\sigma^{-1})^* \tr x^*
    =  (S_\h g)^* \tr (\sigma^{-1} \tr x^*)
    = (S_\h g)^* \tr x^{*_\h} \,,
  \end{split}
  \end{equation}
  where we have used the properties of $\sigma$ from
  Prop.~\ref{th:sigmaunique}.
\end{proof}

\subsection{Representations}
\label{sec:Representations}

The $\h$-adic representations of the algebra $\Ug[[\h]]$ we are most
interested in are those which are the $\h$-adic completion of the
representations of $\Ug$. If $D = (V, \rho)$ is a $\Ug$-module with
complex vector space $V$ and structure homomorphism $\rho: \Ug
\rightarrow \mathrm{End}_{\mathbb{C}}(V)$, its $\h$-adic completion
$\bar{D}=(\bar{V},\bar{\rho})$ is defined on $\bar{V} = V[[\h]]$ with
an order by order extension of $\rho$, $\bar{\rho}(\sum_k \h^k g_k) :=
\sum_k \h^k \rho(g_k) \in \mathrm{End}_{\mathbb{C}}(V)[[\h]] =
\mathrm{End}_{\mathbb{C}[[\h]]}(V[[\h]])$. In particular, $\bar{V}$ is
free over $\mathbb{C}[[\h]]$ and $\bar{\rho}$ is
$\mathbb{C}[[\h]]$-linear. Note that, even if $D$ is an irreducible
representation of $\Ug$ this is no longer true for $\bar{D}$. For
example, $\h \bar{V}$ would be an invariant subspace of $\bar{V}$.
$\bar{D}$ is irreducible only in the sense that there is no subspace
$U \subset V$ such that $\bar{U} := U[[\h]]$ is an invariant subspace
of $\bar{V}$.

Let $\{E_k, H_k, F_k \,\,|\,\, k = 1,\ldots,n\}$ be a Cartan-Weyl
basis of the semisimple Lie algebra $\mathfrak{g}$, and $\mathfrak{h}$
be the Cartan subalgebra which is generated by $\{H_k\}$. Clearly,
$\mathfrak{h}$ is a Lie subalgebra of the enveloping algebra,
$\mathfrak{h} \subset \Ug[[\h]]$.  By construction, the $\h$-adic
Drinfeld-Jimbo deformation of $\Ug$ does contain the Cartan subalgebra
as a Lie subalgebra, as well, $\mathfrak{h} \subset \Uhg$. Drinfeld
has shown, that the isomorphism of $\Ug[[\h]]$ and $\Uhg$ can be
chosen to leave the Cartan subalgebra invariant:

\begin{Theorem}[Drinfeld \cite{Drinfeld:1989}, Prop.~4.3]
  \label{th:CartanPreserve}
  Let $\mathfrak{g}$ be a semisimple Lie algebra and $\mathfrak{h}
  \subset \mathfrak{g}$ a Cartan subalgebra. Then there exists an
  isomorphism of $\h$-adic algebras $\alpha:\Uhg \rightarrow
  \Ug[[\h]]$ such that $\alpha = \id + \Ocal(\h)$ and $\alpha
  \rvert_{\mathfrak{h}} = \id_{\mathfrak{h}}$.
\end{Theorem}

This has important consequences for the representation theory. Recall,
that every finite dimensional irreducible representation of $\Ug$ is a
highest weight representation $D^j = (V^j, \rho^j)$, generated by a
highest weight vector $v \in V^j$ with $\rho^j(E_k)v = 0$ and
$\rho^j(H_k)v = j_k v$ for all $k$, where $j_k \in
\tfrac{1}{2}\mathbb{N}_0$, $j := (j_1,\ldots,j_n)$ being called the
highest weight. Furthermore, there is a basis of $V^j$ which consists
of simultaneous eigenvectors of $H_k$, called the weight basis. The
same is true for the finite dimensional representations of the
Drinfeld-Jimbo algebra $\Uqg$ for a fixed value of $q$
\cite{Lusztig:1988,Rosso:1988}, the weight-$j$ representation $D^j_q =
(V^j, \rho^j_q)$ of $\Uqg$ being defined on the same weight basis as
$D^j$. By the substitution $q \mapsto \E^\h$, $D^j_q$ can be extended
to an $\h$-adic representation $D_\h^j := (V^j[[\h]], \rho_\h^j)$.
Since, $\Uhg$ and $\Ug[[\h]]$ are isomorphic as algebras, there is a
bijection between their representations.
Theorem~\ref{th:CartanPreserve} implies that the isomorphism $\alpha:
\Uhg \rightarrow \Ug[[\h]]$ can be chosen, such that
\begin{equation}
\label{eq:alphaRep}
  \rho_\h^j = \bar{\rho}^j \circ \alpha  \,.
\end{equation}
If both $\rho_\h^j$ and $\bar{\rho}^j$ are known this equation can be
used to calculate the isomorphism $\alpha$.

\section{Constructing Covariant Star Products}
\label{sec:StarProducts}

\subsection{The General Approach}

As explained in Sec.~\ref{sec:Quest} we are asking, if there are
Drinfeld twists which implement the product of quantum spaces by
Eq.~\ref{eq:StarProd2}. To our knowledge, no Drinfeld twist for the
Drinfeld-Jimbo quantum enveloping algebra of a semisimple Lie algebra
has ever been computed. This suggests, that it will be rather
difficult to answer this question on an algebraic level. The
representations of Drinfeld twists, however, can be computed as we
will demonstrate for $\Uhsu$ in Prop.~\ref{th:Frep1}. Therefore, we
propose the following approach, which allows us to tackle the problem
on a representation theoretic level:

Consider a $\Uhg$-covariant quantum space algebra $\Xcal_h$ and its
undeformed limit, the $\Ug$-covariant space algebra $\Xcal$.

\begin{enumerate}
  
\item Determine the irreducible highest weight representations of all
  possible Drinfeld twists from $\Delta$ to $\Delta_\h$. 
  
\item Determine the basis $\{ T_{m,k}^j \}$ of the quantum space
  $\Xcal_\h$ which completely reduces $\Xcal_\h$ into irreducible
  highest weight representations of $\Uhg$,
  \begin{equation}
    \Xcal_\h \cong \bigoplus_{j,k}
    \operatorname{Span}_{\mathbb{C}[[\h]]}
    \{ T^j_{m,k} \,\,| \,\,m \text{ weight of } D^j_\h\} \,,
  \end{equation} 
  such that $g \tr T^j_{m,k} = T^j_{m',k} \rho_\h^j(g)^{m'}{}_{m}$ for
  all $g \in \Uhg$, where $m$ is a weight, $j$ is the highest weight,
  $\rho^j_\h$ is the structure map of the $\h$-adic highest weight-$j$
  representation $D^j_\h$ of $\Uhg$ as explained in
  Sec.~\ref{sec:Representations}, and where $k$ labels the possibly
  degenerate highest weight-$j$ subrepresentations.
  
\item Calculate the multiplication map $\mu_\h$ of $\Xcal_\h$ with
  respect to this basis. The undeformed limit $T^{\prime j}_{m,k} :=
  \lim_{\h \rightarrow 0} T^j_{m,k}$ then yields the basis which
  completely reduces the undeformed space algebra $\Xcal$. The limit
  $\mu = \lim_{\h \rightarrow 0} \mu_\h$ is the commutative
  multiplication map with respect to this basis.
  
\item With respect to the basis $\{ T^{\prime j_1}_{m_1,k_1} \otimes
  T^{\prime j_2}_{m_2,k_2} \}$ of $(\Xcal \otimes \Xcal)[[\h]]$ the
  action of the twist is given by the highest weight representations
  $(\rho^{j_1}\otimes \rho^{j_2})(\F)$. Now we can check if one of the
  twists realizes the deformed multiplication map by
  Eq.~\eqref{eq:StarProd2} as linear operator with respect to this
  basis.

\end{enumerate}
Since this procedure reduces the algebraic problem to a representation
theoretic one, it works well for quantum spaces of $\Uhsu$, $\Uhso$,
and $\UhslC$ where the representation theory is well understood.


\subsection{The Drinfeld Twists of $\Uhsu$}

We now consider the case of $\mathfrak{g} = A_1$, the complex Lie
algebra with Cartan-Weyl basis $\{E,H,F\}$ and relations $[H,E]=2E$,
$[H,F]=-2F$, $[E,F] = H$. The real form of $A_1$ which corresponds to
the $*$-structure $E^* = F$, $H^* = H$, $F^* = E$ is $\mathrm{su_2}$,
the Lie algebra of the group of unitary 2$\times$2-matrices.

\begin{Definition}
  The complex $\h$-adic algebra generated by $E$, $H$, $F$ with
  commutation relations
  \begin{xalignat}{3}
    [H,E] &= 2E \,,& [H,F] &= -2F \,, &
    [E,F] &= \frac{\E^{\h H}-\E^{-\h H}}{\E^\h-\E^{-\h}} \,,
  \end{xalignat}
  Hopf structure
  \begin{xalignat}{3}  
  \Delta'(E)&= E\otimes \E^{\h H} + 1\otimes E \,,&
  S'(E)&=-E\E^{-\h H} \,,& \varepsilon'(E)&=0 \notag\\
  \Delta'(F) &= F\otimes 1 + \E^{-\h H} \otimes F \,, &
  S'(F)&=-\E^{\h H}F \,,& \varepsilon'(F)&=0 \,, \\ 
  \Delta'(H)&= H\otimes 1 + 1\otimes H \,,&
  S'(H)&= -H \,,& \varepsilon'(H)&=0 \,, \notag
  \end{xalignat}
  and involution $E^{*'} = F \E^{\h H}$, $F^{*'} = \E^{-\h H}E$ ,
  $H^{*'} = H$ is called $\Uhsu$, the $\h$-deformation of $\Usu$
  \cite{Kulish:1983,Sklyanin:1985}.
It is quasitriangular with universal
$\R$-matrix \cite{Drinfeld:1986}
\begin{equation}
\label{eq:Rmatrix}
  \R = \E^{\h (H\otimes H)/2} \sum_{n=0}^{\infty}
  \E^{\h n(n-1)/2} \frac{(\E^\h- \E^{-\h})^n}{[n]!} (E^n \otimes F^n) \,.
\end{equation}
\end{Definition}

By construction, the commutation relations and Hopf $*$-structure maps
of $\Uhsu$ coincide in zeroth order of $\h$ with those of $\Usu$.
Therefore, $\Uhsu$ is a deformation of $\Usu$ as Hopf $*$-algebra. The
$\h$-adic deformation is obtained from the $q$-deformation by the
substitutions $q = \E^\h$, $K = \E^{\h H}$, and $K^{-1} = \E^{-\h H}$.
By the same substitution we obtain for each
$j\in\tfrac{1}{2}\mathbb{N}_0$ the $h$-adic spin-$j$
$*$-representation
\begin{equation}
\label{eq:Irrepsh}
\begin{aligned}
  \rho^j_\h(E)\Ket{j,m} &= \E^{\h(m+1)} \sqrt{[j+m+1][j-m]}
                           \,\Ket{j,m+1} \\
  \rho^j_\h(F)\Ket{j,m} &= \E^{-\h m} \sqrt{[j+m][j-m+1]}\,\Ket{j,m-1}\\
  \rho^j_\h(H)\Ket{j,m} &= 2m\Ket{j,m} \,, 
\end{aligned}
\end{equation}
on the $(2j+1)$-dimensional free $\mathbb{C}[[\h]]$-module $V^j[[\h]]$
with orthonormal weight basis $\{\Ket{j,m}, m= -j,-j+1,\ldots,j\}$,
which we denote by $D_\h^j := (V^j[[\h]],\rho_\h^j)$. Using the
coproduct, tensor representations are constructed as
\begin{equation}
  D_\h^{j_1} \otimes D_\h^{j_2} := \bigl( V^{j_1} \otimes V^{j_2},
  \rho_\h^{j_1 \otimes j_2}
  := (\rho_\h^{j_1} \otimes \rho_\h^{j_2})\circ \Delta' \bigr)
\end{equation}
and analogously for the undeformed case. The decomposition of such a
tensor representation into its irreducible subrepresentations is the
Clebsch-Gordan series
\begin{equation}
\label{eq:CGSeries1}
  D_\h^{j_1}\otimes D_\h^{j_2} \cong
  D_\h^{|j_1-j_2|} \oplus D_\h^{|j_1-j_2|+1}
  \oplus \ldots \oplus D_\h^{j_1+j_2} \,.
\end{equation}
Let us denote the embedding of the irreducible spin-$j$ component into
the tensor representation by $C_q^{j_1 j_2 j}$ and the projection onto
this component by $(C_q^{j_1 j_2 j})^{-1}$, such that
\begin{equation}
\label{eq:CGSeries3}
  \rho_\h^{j}(g) \,C_q^{j_1 j_2 j} = C_q^{j_1 j_2 j} \,
  \rho_\h^{j_1 \otimes j_2}(g)
  = C_q^{j_1 j_2 j} \,
  (\rho_\h^{j_1} \otimes \rho_\h^{j_2})(\Delta' g )
\end{equation}
for all $g \in \Uhsu$. Denoting the basis vectors of
$D_\h^{j_1}\otimes D_\h^{j_2}$ by $\Ket{j_1,m_1;j_2,m_2}$ and those of
the irreducible spin-$j$ subrepresentation by $\Ket{j,m}$, the
$q$-Clebsch-Gordan coefficients are defined as
\begin{equation}
\label{eq:CGSeries4}
  \CGqs{j_1}{j_2}{j}{m_1}{m_2}{m} :=
  \Bra{j_1,m_1;j_2,m_2}  C_q^{j_1 j_2 j} \Ket{j,m} \,. 
\end{equation}
The $q$-Clebsch-Gordan coefficients are not unique, because the basis
vectors $\Ket{j,m}$ are only determined up to a phase. We will follow
the choice of \cite{Schmuedgen}, where
\begin{equation}
\label{eq:CGSeries5}
\begin{aligned}
  \sum_{m_1, m_2} \CGqs{j_1}{j_2}{j}{m_1}{m_2}{m}
  \CGqs{j_1}{j_2}{j'}{m_1}{m_2}{m'} &= \delta_{mm'} \delta_{jj'} \\
  \sum_{j, m} \CGqs{j_1}{j_2}{j}{m_1}{m_2}{m}
  \CGqs{j_1}{j_2}{j}{m_1'}{m_2'}{m} &= \delta_{m_1m_1'} \delta_{m_2
    m_2'}\,.
\end{aligned}
\end{equation}
Clearly, the representation theory of $\Uhsu$ is a deformation of the
one of $\Usu$. In the limit $\h \rightarrow 0$ or, equivalently, $q
\rightarrow 1$ of Eq.~\eqref{eq:CGSeries4} we get back the undeformed
Clebsch-Gordan coefficients. Since $\Uhsu$ and $\Usu[[\h]]$ are
isomorphic as $*$-algebras their representations are isomorphic, as
well. Due to Theorem~\ref{th:CartanPreserve} we can choose the
isomorphism $\alpha: \Uhsu \rightarrow \Usu[[\h]]$ as in
Eq.~\eqref{eq:alphaRep} such that
\begin{equation}
\label{eq:alphastar}
  \rho^j_\h = \rho^j \circ \alpha \,,
\end{equation}
where $\rho^j: \Usu[[\h]] \rightarrow \mathrm{End}(V^j)[[\h]]$ is the
$\h$-adically extended structure map of the undeformed enveloping
algebra (We omit the bar which denoted the $\h$-adic completion in
Eq.~\eqref{eq:alphaRep}). If we set $g = \alpha^{-1}(g')$ in
Eq.~\eqref{eq:CGSeries3} we thus get
\begin{equation}
\label{eq:CGSeries6}
  \rho^{j}(g') \,C_q^{j_1 j_2 j} = C_q^{j_1 j_2 j} \,
  (\rho^{j_1} \otimes \rho^{j_2})(\Delta_\h\, g' )
\end{equation}
for all $g' \in \Usu$, where $\Delta_\h$ is defined as in
Eq.~\eqref{eq:HopfDeform}.

It is rather obvious that the representations of Drinfeld twists
should be given by a contraction of the deformed and undeformed
Clebsch-Gordan coefficients, as it was already mentioned in
\cite{Curtright:1991}.
\begin{Proposition}
  \label{th:Frep1}
  Let $\F$ be a counital Drinfeld twist from $\Usu$ to $\Uhsu$. The
  irreducible representations of $\F$ are of the form
  \begin{equation}
  \label{eq:Frep1}  
    (\rho^{j_1}\otimes \rho^{j_2})(\F)^{m_1 m_2}{}_{m_1' m_2'}
    = \sum_{j,m}\eta(j_1,j_2,j)
    \CGqs{j_1}{j_2}{j}{m_1}{m_2}{m}
    \CGs{j_1}{j_2}{j}{m_1'}{m_2'}{m} \,,
  \end{equation}
  where for given values of $j_1$, $j_2$, and $j$ the factor
  $\eta(j_1,j_2,j) \in \mathbb{C}[[\h]]$ is a formal power series in
  $\h$ with $\eta(j_1,j_2,j) = 1 + \Ocal(\h)$.
\end{Proposition}
\begin{proof}
  Within a $D^{j_1} \otimes D^{j_2}$ tensor representation we get for
  all $g\in \Usl$
\begin{equation}
\label{eq:Frep1b}
\begin{split}
  C_q^{j_1j_2 j}\,&(\rho^{j_1}\otimes \rho^{j_2})(\F)\,
  (C^{j_1j_2 j'})^{-1} \rho^{j'}(g) \\
  &=
  C_q^{j_1j_2 j}\,(\rho^{j_1}\otimes \rho^{j_2})(\F \Delta g)
  (C^{j_1j_2 j'})^{-1} \\
  &=
  C_q^{j_1j_2 j}\,(\rho^{j_1}\otimes \rho^{j_2})
  ((\Delta_\h\, g) \F)(C^{j_1j_2 j'})^{-1} \\
  &=
  \rho^{j}(g)C_q^{j_1j_2 j}\,
  (\rho^{j_1}\otimes \rho^{j_2})(\F)(C^{j_1j_2 j'})^{-1} \,,
\end{split}
\end{equation}
where we have used Eq.~\eqref{eq:CGSeries6} and the analogous relation
for the undeformed case. Let us develop $\eta := C_q^{j_1j_2
  j}\,(\rho^{j_1}\otimes \rho^{j_2})(\F)\, (C^{j_1j_2,j'})^{-1}$ into
an $\h$-adic series $\eta = \sum_k \h^k \eta_k$. Then
Eq.~\eqref{eq:Frep1b} implies that each $\eta_k$ is a module map from
the spin-$j'$ to the spin-$j$ irreducible subrepresentation of the
$D^{j_1} \otimes D^{j_2}$ tensor representation. By Schur's lemma each
$\eta_k$ must be zero for $j\neq j'$, while for $j = j'$ the $\eta_k$
are $\mathbb{C}[[\h]]$-scalar multiples of the identity map
$\id_{D^j}$. Hence,
\begin{equation}
\begin{aligned}
  C_q^{j_1j_2,j}\,(\rho^{j_1}\otimes \rho^{j_2})(&\F)\,
  (C^{j_1j_2,j'})^{-1} =  \eta(j_1,j_2,j) \, \delta_{jj'}
  \,\id_{D^j} \\ \Leftrightarrow \quad
  (\rho^{j_1}\otimes \rho^{j_2})(&\F)\, =
  \sum_j \eta(j_1,j_2,j) (C_q^{j_1j_2,j})^{-1} C^{j_1j_2,j} \,,
\end{aligned}
\end{equation}
where for given $j_1$, $j_2$, $j$, $\eta(j_1,j_2,j) \in
\mathbb{C}[[\h]]$. Taking matrix elements of the last equation and
using the definition~\eqref{eq:CGSeries4} of the Clebsch-Gordan
coefficients yields Eq.~\eqref{eq:Frep1}. Finally, $\F = 1 +
\Ocal(\h)$ implies $\eta(j_1,j_2,j) = 1 + \Ocal(\h)$.
\end{proof}

The Drinfeld twist with the canonically simplest representations would
be the one with $\eta(j_1,j_2,j) = 1$. Such a twist exists, indeed.

\begin{Proposition}
\label{th:Frep2}
There is a unique Drinfeld twist $\Fs$, called the standard twist, for
which the irreducible representations of Prop.~\ref{th:Frep1} are such
that $\eta(j_1,j_2,j)=1$ for all $j_1$, $j_2$, $j$.
\end{Proposition}

\begin{proof}
  Define a scalar conjugation $u \mapsto \overline{u}$ on $\Usu$ by
  extending the identity map on the Cartan-Weyl generators of
  $\mathrm{su}_2$ to a conjugate linear automorphism of $\Usu$, e.g.,
  $\overline{\alpha EF} = \overline{\alpha}EF$ etc. Since both
  coproducts are real with respect to this conjugation,
  $\overline{\Delta(g)} = \Delta(\overline{g})$ and
  $\overline{\Delta_\h(g)} = \Delta_\h(\overline{g})$ the conjugation
  of Eq.~\eqref{eq:StarProd4} shows that, if $\F$ is a counital
  Drinfeld twist, so is $\overline{\F}$ and, hence, $\F':=
  \frac{1}{2}(\F + \overline{\F})$. In the
  representation~\eqref{eq:Irrepsh} and its undeformed limit the
  Cartan-Weyl generators are represented by real matrices, thus,
  $\rho^j(\overline{g})^m{}_{m'} = \overline{\rho^j(g)^m{}_{m'}}$. We
  conclude that if $\eta(j_1,j_2,j) = \eta$ are the factors of the
  representations of $\F$, $\overline{\eta}$ are those of
  $\overline{\F}$ and $\eta'= \eta + \overline{\eta}$ those of $\F'$.
  As in the proof of Prop.~\ref{th:Forth} the twist $\F'' :=
  \F'(F^{\prime -1}(F^{\prime *})^{-1})^{1/2}$ is unitary. The factors
  $\eta''$ of the representations of $\Fs$ are real because those of
  $\F'$ are real, and unitary $\overline{\eta''} = \eta^{\prime\prime
    -1}$ because $\F''$ is unitary. Hence, $\eta = \eta(j_1,j_2,j) =
  1$.  Since any two $\F$ with the same representations are equal,
  $\F''$ is the unique standard twist.
\end{proof}

\begin{Proposition}
  The standard twist $\Fs$ is orthogonal.
\end{Proposition}

\begin{proof}
\label{th:Frep3}
  It was shown in the proof of Prop.~\ref{th:Frep2} that $\Fs$ is
  unitary. Moreover,
  \begin{multline}
  \label{eq:Frep3}
    (\rho^{j_1}\otimes \rho^{j_2})
    \bigl((S \otimes S)(\Fs) \bigr)^{m_1 m_2}{}_{m_1' m_2'}
    = \sum_{j,m} \CGqs{j_1}{j_2}{j}{-m_1'}{-m_2'}{m} \,
      \CGs{j_1}{j_2}{j}{-m_1}{-m_2}{m} \\
    = \sum_{j,m} \CGs{j_2}{j_1}{j}{m_2}{m_1}{m} \,
    \CGqs{j_2}{j_1}{j}{m_2'}{m_1'}{m}
    = (\rho^{j_1}\otimes \rho^{j_2})(\Fs_{21}^{-1})^{m_1
      m_2}{}_{m_1' m_2'} \,,
  \end{multline}
  where we have used $\rho^j(Sg)^{m}{}_{m'} = (-1)^{m-m'}
  \rho^j(g)^{-m'}{}_{\!\!-m}$ and that
  $\CGqs{j_1}{j_2}{j}{m_1}{m_2}{m} =
  \CGqs{j_2}{j_1}{j}{-m_2}{-m_1}{-m}$. We conclude that $(S \otimes
  S)(\Fs) = \Fs_{21}^{-1} = (* \otimes *)(\Fs_{21})$.
\end{proof}

\subsection{The Quantum Plane}

The quantum plane \cite{Manin:1988} is perhaps the simplest nontrivial
example of a homogeneous quantum space. In analogy to the undeformed
case, the generators $x_-$ and $x_+$ are defined to carry the
fundamental spin-$\frac{1}{2}$ representation of $\Uhsu$,
\begin{equation}
  g \tr x_m := x_{m'}\rho_\h^{\frac{1}{2}}(g)^{m'}{}_{m} \,,
\end{equation} 
where the indices run through $\{-,+\} = \{-\frac{1}{2},
+\frac{1}{2}\}$ (summation over repeated upper and lower indices). We
also denote the generators by $x \equiv x_-$ and $y \equiv x_+$. Let
$\mathbb{C}\langle x_-, x_+ \rangle[[\h]]$ be the free $\h$-adic
algebra generated by $x_-$ and $x_+$. By construction, an algebra
which is freely generated by a $\Uhsu$-module is a $\Uhsu$-module
algebra. The quadratic terms $x_{m_1}x_{m_2}$ thus carry a
spin-$(\frac{1}{2} \otimes \frac{1}{2})$ tensor representation. If we
want to divide the free algebra by quadratic relations in such a way
that the quotient algebra is again a $\Uhsu$-module algebra we must
divide by an ideal which is generated by a submodule of the
representation $D_\h^{1/2} \otimes D_\h^{1/2} \cong D_\h^{0} \oplus
D_\h^{1}$ of all quadratic terms. Dividing by $D_\h^1$ would yield a
deformation of the exterior algebra, whereas dividing by the scalar
part $D_\h^0$ yields the desired deformation of the commutative
algebra of functions on the 2-dimensional plane. This amounts to the
commutation relations
\begin{equation}
\label{eq:qPlaneCommute}
  \sum_{m_1, m_2} \CGqs{1/2}{1/2}{0}{m_1}{m_2}{0}\, x_{m_1} x_{m_2} = 0
  \quad\Leftrightarrow\quad x y = q y x \,,
\end{equation}
where $q = \E^\h$. 

\begin{Definition}
  The $\h$-adic algebra freely generated by $x \equiv x_-$ and $y
  \equiv x_+$ with commutation relations~\eqref{eq:qPlaneCommute} is
  called the $\h$-adic quantum plane $\Xcal_\h(\mathbb{C}^2)$.
\end{Definition} 

We now want to write the product of the quantum plane, $\mu_\h(x_1
\otimes x_2) := x_1 x_2$, explicitly as a linear map with respect to a
basis. For our purposes, the appropriate choice is the basis which
reduces the quantum plane as $\Uhsu$-module into its irreducible
subrepresentations. In order to find such a basis, we recall that, as
in the undeformed case, finding the irreducible spin-$j$
subrepresentations is the matter of finding the highest weight-$j$
vectors, that is, the elements of $\Xcal_\h(\mathbb{C}^2)$ which
transform as $\Ket{j,j}$ in Eqs.~\eqref{eq:Irrepsh}. A simple ansatz
shows that, up to scalar multiples, the only element of the quantum
plane with this property is $x_+^{2j}$. Acting on $x_+^{2j}$ with the
ladder operator $F$ generates the other basis vectors of the
$D^j_\h$-subrepresentation. Identifying in Eqs.~\eqref{eq:Irrepsh}
$\Ket{j,m}$ with $T^j_m$, we have to define
\begin{equation}
\label{eq:TensorBasis1}
  T^j_m := q^{\frac{1}{2}(j-m)(2m-j+1)}
  \sqrt{\frac{[j+m]!}{[2j]![j-m]!}}\, (F^{j-m} \tr x_+^{2j})
\end{equation}
for $m \in \{ -j, -j+1, \ldots, j\}$, such that
\begin{equation}
  g \tr T^j_m = T^j_{m'}\rho_\h^j(g)^{m'}{}_m
\end{equation}
for all $g \in \Uhsu$. The basis $\{T^j_m\}$ then reduces the quantum
plane into its irreducible subrepresentations, $\Xcal_\h(\mathbb{C}^2)
= D^{0}_\h \oplus D^{1/2}_\h \oplus D^1_\h \ldots$ Calculating
Eq.~\eqref{eq:TensorBasis1} explicitly, yields
\begin{equation}
\label{eq:TensorBasis2}
  T^j_m =
  \qBinom{2j}{j+m}_{q^{-2}}^{\frac{1}{2}} x_-^{j-m} x_+^{j+m}
  \,,\quad\text{where}\quad
  \qBinom{n}{k}_{q^{-2}} :=
   q^{k(k-n)}\,\frac{[n]!}{[n-k]![k]!}
\end{equation}
is the $q$-binomial coefficient. By construction, $T^{j_1}_{m_1}
T^{j_2}_{m_2}$ carries a spin-$(j_1 \otimes j_2)$ tensor
representation which can be reduced using the $q$-Clebsch-Gordan
coefficients. Hence, the elements
\begin{equation}
\label{eq:Adef1}
  A^j_m := \sum_{m_1,m_2}\CGqs{j_1}{j_2}{j}{m_1}{m_2}{m}
    \, T^{j_1}_{m_1} T^{j_2}_{m_2}
\end{equation} 
are either zero or the basis of a spin-$j$ subrepresentation of
$\Xcal_\h(\mathbb{C}^2)$. Since $T^j_m$ generates the only spin-$j$
subrepresentation, $A^j_m$ must be proportional to $T^j_m$. Moreover,
because of its homogeneous commutation relations the algebra
$\Xcal_\h(\mathbb{C}^2)$ is graded. That is, the degree of the product
of two homogeneous elements $T^{j_1}_{m_1}$ and $T^{j_2}_{m_2}$ is the
sum of their degrees, $\mathrm{deg} (T^{j_1}_{m_1} T^{j_2}_{m_2}) =
\mathrm{deg}(T^{j_1}_{m_1}) + \mathrm{deg}(T^{j_2}_{m_2}) = 2(j_1 +
j_2)$. As $A^j_m$ and $T^j_m$ are proportional, they must have the
same degree. Thus, $A^j_m$ has to vanish unless $j = j_1 + j_2$.
Looking at the highest weight vectors we find $A^{j_1+j_2}_{j_1+j_2} =
x_+^{2(j_1 + j_2)} = T^{j_1 + j_2}_{j_1+j_2}$. Using the orthogonality
relation~\eqref{eq:CGSeries5} we can move the $q$-Clebsch-Gordon
coefficients to the left hand side of Eq.~\eqref{eq:Adef1}. As end
result, the multiplication map $\mu_\h$ of the quantum plane is given
with respect to the basis $\{T^j_m\}$ by
\begin{equation}
\label{eq:qPlaneMulti}
  \mu_\h(T^{j_1}_{m_1} \otimes T^{j_1}_{m_1}) =
  \CGqs{j_1}{j_2}{j_1 + j_2}{m_1}{m_2}{m_1 + m_2}
  \,\,T^{j_1+j_2}_{m_1+m_2} \,.
\end{equation}
For the undeformed limit of the basis $T^{\prime j}_m := \lim_{\h
  \rightarrow 0} T^j_m$ and multiplication map $\mu := \lim_{\h
  \rightarrow 0} \mu_\h$ one gets
\begin{equation}
\label{eq:PlaneUndef}
  T^{\prime j}_m =
  \tbinom{2j}{j+m}^{\frac{1}{2}} x_-^{j-m} x_+^{j+m}
  \,,\quad
  \mu(T^{\prime j_1}_{m_1} \otimes T^{\prime j_1}_{m_1}) =
  \CGs{j_1}{j_2}{j_1 + j_2}{m_1}{m_2}{m_1 + m_2}
  \,\,T^{\prime j_1+j_2}_{m_1+m_2}
\end{equation}
Comparing Eqs.~\eqref{eq:qPlaneMulti} and \eqref{eq:PlaneUndef} with
the representations~\eqref{eq:Frep1} of the Drinfeld twists we obtain
the following

\begin{Proposition}
\label{th:PlaneF}
  Let $\mu_\h$ be the multiplication map~\eqref{eq:qPlaneMulti} of the
  $\h$-adic quantum plane $\Xcal_\h(\mathbb{C}^2)$, $\mu = \lim_{\h
    \rightarrow 0} \mu_\h$ its undeformed limit, and $\Fs$ the
  standard twist of Prop.~\ref{th:Frep2}. Then $\mu_\h$ is the
  deformation~\eqref{eq:StarProd2} of $\mu$ by $\Fs$, $\mu_\h(x
  \otimes y) = \mu(\Fs^{-1} \tr [x \otimes y])$.
\end{Proposition}

Finally, let us turn to real structures. Since according to
Prop.~\ref{th:Frep3} $\Fs$ is real in the sense of
Eq.~\eqref{eq:Freal}, Prop.~\ref{th:StarModuleDeform} applies. Within
$\Uhsu$ we have $S^{\prime 2}(g) = K g K^{-1}$ with $K = \E^{\h H}$
for all $g \in \Uhsu$. Clearly, $K$ is group-like, $\Delta'(K) = K
\otimes K$. Recall that, in order for the representations of
$\Usu[[\h]]$ and $\Uhsu$ to be related by Eq.~\eqref{eq:alphastar}, we
chose the isomorphism $\alpha:\Uhsu \rightarrow \Usu[[\h]]$ according
to Theorem~\ref{th:CartanPreserve} such that $\alpha(H) = H$. On the
one hand, $S_\h^2(g) = (\alpha \circ S^{\prime 2} \circ \alpha)(g) = K
g K^{-1}$. On the other hand, $S_\h^2(g) = \sigma^{2}g\sigma^{-2}$,
for the unique element $\sigma$ from Prop.~\ref{th:sigmaunique}.
Hence, $K \sigma^{-2} = \sigma^{-2} K$ is central, so $S_\h(g) =
K^{1/2} (Sg) K^{-1/2}$. Since $K^{1/2} = 1 + \Ocal(\h)$ is group-like
the uniqueness of $\sigma$ implies that $\sigma = K^{1/2} = \E^{\h
  H/2}$. Prop.~\ref{th:StarModuleDeform} now tells us that for a given
covariant real structure $*$ on the undeformed space algebra $\Xcal
\equiv \Xcal(\mathbb{C}^2) = \mathbb{C}[x_-,x_+][[\h]]$ of the plane,
we have to define the deformed $*$-structure $*_\h = * + \Ocal(\h)$ by
\begin{equation}
\label{eq:PlaneStar}
  x^{*_\h} = \E^{- \frac{\h H}{2}} \tr x^*
\end{equation}
such that $(\Xcal, \mu_\h, *_\h)$ becomes a module $*$-algebra of
$\Usu[[\h]]$ with respect to the deformed Hopf structure.

\subsection{The Quantum Lorentz Algebra}

We recall the definition of quantum Euclidean algebra in 4 dimensions
and the quantum Lorentz algebra.

\begin{Definition}
\label{th:LorentzAlgebra}
The tensor product Hopf $*$-algebra $\Uhsu \otimes \Uhsu$ is the
$\h$-adic quantum enveloping algebra of $\mathrm{so_4}$, $\Uhso$. Let
$\R$ be the universal $\R$-matrix~\eqref{eq:Rmatrix} of $\Uhsu$. The
Hopf algebra obtained by twisting $\Uhso$ with $\R^{-1}_{23} = 1
\otimes \R^{-1} \otimes 1$ together with the $*$-structure
\begin{equation}
\label{eq:LorentzStar}
  (a\otimes b)^* =\R_{21}(b^*\otimes a^*)\R_{21}^{-1}
\end{equation}
for $a, b \in \Uhsu$ is the $\h$-adic quantum Lorentz algebra $\UhslC$
\cite{Majid:1993}.
\end{Definition}

The tensor Hopf $*$-structure of $\Uhso$ is given by
$\varepsilon^{\otimes 2} := \varepsilon \otimes \varepsilon$,
$S^{\otimes 2} := S \otimes S$, $\Delta^{\!\otimes 2} := \tau_{23}
\circ (\Delta \otimes \Delta)$, where $\tau$ is the flip of the tensor
factors, $\tau(a \otimes b) = b \otimes a$, and $*^{\otimes 2} := *
\otimes *$. Looking at Cor.~\ref{th:HopfTwist}, (v) for the twist $\F
:= \R^{-1}_{23}$ we find that the coassociator is the unit, $\Phi =
1$. Twists with unital coassociator are 2-cocycles in the sense of
\cite{Majid}. The cocycle property guarantees that the twisted
coproduct is coassociative. For the antipode $S(a \otimes b) =
\sigma_1 (Sa\otimes Sb) \sigma_1^{-1}$ we have to compute
$\sigma_1^{-1} = S(\F_{[1]}^{-1})\F_{[2]}^{-1} = (1 \otimes
S\R_{[1]})(\R_{[2]} \otimes 1) = \R^{-1}_{21}$. Whence, the Hopf
structure of $\UhslC$ reads explicitly
\begin{equation}
\label{eq:LorentzHopf}
  \Delta(a\otimes b) =\R^{-1}_{23}
    \Delta^{\!\otimes 2} (a\otimes b)\R_{23} \,,\quad
    S(a\otimes b) = \R_{21} (Sa\otimes Sb) \R^{-1}_{21} \,,
\end{equation}
while according to Cor.~\ref{th:HopfTwist}, (ii) the counit stays
undeformed.

In the undeformed case $\mathrm{so_4}$ and $\slC$ are real forms of
the same complex Lie algebra $A_1 \otimes A_1$, that is, $\Uso$ and
$\UslC$ differ only by their $*$-structure, whereas in the
$q$-deformed case the Hopf structures differ, as well. The reason for
introducing the twist in the Hopf structure of $\UhslC$ is that only
then the quantum Lorentz algebra contains a Hopf $*$-subalgebra of
rotations, embedded by the coproduct $\Delta:\Uhsu \hookrightarrow
\UhslC$, which is an essential feature for its physical
interpretation. The $*$-structure of $\UslC$ is not a twist of the
product $*$-structure of $\Uhso$, but of the flipped $*$-structure
$\tau \circ (* \otimes *)$. While Eq.~\eqref{eq:LorentzStar} clearly
defines an algebra antihomomorphism, the involution property $*^2 =
\id$ relies on the additional property $\R^{* \otimes *} = \R_{21}$ of
the $\R$-matrix of $\Uhsu$.

The twisting from $\Uhso$ to $\UhslC$ can be extended to Drinfeld
twists and module algebras.

\begin{Corollary}
\label{th:MinkowskiTwist}
  Let $\F'$ and $\F''$ be Drinfeld twists from $\Usu$ to $\Uhsu$, $\R$
  the universal $\R$-matrix of $\Uhsu$, $\Xcal$ an $\h$-adic $\Uhsu
  \otimes \Uhsu$-module.
\begin{itemize}
  
\item[\textup{(i)}] $\F'_{13} \F''_{24}$ is a Drinfeld twist of from
  $\Uso$ to $\Uhso$.
  
\item[\textup{(ii)}] $\R^{-1}_{23}\F'_{13} \F''_{24}$ is a Drinfeld
  twist of from $\UslC$ to $\UhslC$.
  
\item[\textup{(iii)}] If $(\Xcal, \mu)$ is a module algebra of $\Uhso$
  with multiplication map $\mu$, then $(\Xcal, \tilde{\mu})$ with the
  twisted multiplication defined as $\tilde{\mu}(x \otimes y) :=
  \mu(\R_{23} \tr [ x \otimes y] )$ is a module algebra of $\UhslC$.
  
\item[\textup{(iv)}] If $(\Xcal,\mu, *)$ is a module $*$-algebra of
  the Hopf algebra $\Uhso$ with flipped $*$-structure $\tau \circ (*
  \otimes *)$ then $(\Xcal, \tilde{\mu}, *)$ is a module $*$-algebra
  of $\UslC$.

\end{itemize}
\end{Corollary}

\begin{proof}
(i) $\F'_{13} \F''_{24} = \tau_{23}(\F' \otimes \F'')$ is clearly a
tensor product twist of the tensor coproduct $\Delta^{\otimes 2} =
\tau_{23}\circ (\Delta \otimes \Delta)$
  
(ii) Twisting in two steps by $\F'_{13} \F''_{24}$ from $\UslC = \Uso$
to $\Uhso$ and then by $\R^{-1}_{23}$ from $\Uhso$ to $\UhslC$ is the
same as twisting at once by $\R^{-1}_{23}\F'_{13} \F''_{24}$ from
$\UslC$ to $\UhslC$.
 
(iii) Note, that $\tilde{\mu}$ is the twisted
multiplication~\eqref{eq:StarProd2} for $\F = \R^{-1}_{23}$. All we
have to check is associativity. The defining properties of a universal
$\R$-matrix imply that the Drinfeld coassociator~\eqref{eq:CoassDef}
is equal to the unit. Hence, Eq.~\eqref{eq:Quest1} is trivially
satisfied.
  
(iv) Let us denote the flipped $*$ structure by $*_\tau := \tau \circ
(* \otimes *)$. We verify that $\F = \R^{-1}_{23}$ is real in the
sense of Eq.~\eqref{eq:Freal},
\begin{equation}
\begin{split}
  (S^{\otimes 2} \otimes S^{\otimes 2})(\F) &=
  1 \otimes (S \otimes S)(\R^{-1}) \otimes 1
  = 1 \otimes (\R_{21}^{-1})^{ *\otimes *} \otimes 1 \\ 
  &= (\R_{[2]}^{-1} \otimes 1)^{*_\tau} \otimes
  (1 \otimes \R_{[1]}^{-1})^{*_\tau} 
  = (*_\tau \otimes *_\tau)(\F_{21}) \,,
\end{split}
\end{equation}
where we have used $(S \otimes S)(\R) = \R$ and $\R^{* \otimes *} =
\R_{21}$. Eq.~\eqref{eq:antihom} shows that the $*$-structure of
$\Xcal$ is an antihomomorphism with respect to both multiplications,
$\mu$ and $\tilde{\mu}$. Then we have to check
Eq.~\eqref{eq:StarModule}. Denoting by $S_\slC$ and $*_\slC$ the
antipode and $*$-structure of $\UhslC$, we get
\begin{equation}
\begin{split}
  [S_{\slC}(a \otimes b)]^{*_\slC}
  &= [\R_{21}(Sa\otimes Sb)\R_{21}^{-1}]^{*_\slC} \\ 
  &= \R_{21}[( \R_{21}^{-1}((Sb)^* \otimes (Sa)^*)\R_{21})]\R_{21}^{-1} \\
  &= (Sb)^* \otimes (Sa)^*
  = [S^{\otimes 2}( a \otimes b)]^{*_\tau} \,,
\end{split}
\end{equation}
whence, $[(a \otimes b) \tr x]^* = [S^{\otimes 2}(a \otimes
b)]^{*_\tau} \tr x^* = [S_\slC(a \otimes b)]^{*_\slC} \tr x^*$.
\end{proof}

\subsection{Quantum Minkowski Space}

Quantum Minkowski space $\Minkh$ is a noncommutative deformation of
the function algebra on real world 1+3-dimensional spacetime
\cite{Carow-Watamura:1990}. By definition, $\Minkh$ is the
$\UhslC$-module algebra whose generators carry the fundamental
representation. It was shown in Cor.~\ref{th:MinkowskiTwist} that any
$\UhslC$-module algebra is the twist of a $\Uhso$-module algebra. We
will first compute the multiplication map of this $\Uhso$-module
algebra, quantum Euclidean 4-space, and then twist it to obtain the
multiplication map of quantum Minkowski space.

Because $\Uhso$ is the product Hopf algebra of two $\Uhsu$, any
irreducible representation is the product of two irreducible
representations of $\Uhsu$,
\begin{equation}
  D_\h^{(j,j')} := (V^j \otimes V^{j'},
  \rho_\h^j \otimes \rho_\h^{j'}) \,.
\end{equation}
The generators $X_{mm'}$ of the quantum Euclidean 4-space are defined
to carry the fundamental spin-$(\frac{1}{2},\frac{1}{2})$ representation,
\begin{equation}
\label{eq:XModule}
  (g\otimes g') \tr X_{mm'} = X_{\tilde{m}\tilde{m}'} \,
  \rho_\h^{\frac{1}{2}}(g)^{\tilde{m}}{\!}_{m} \,
  \rho_\h^{\frac{1}{2}}(g')^{\tilde{m}'}{\!}_{m'} \,,  
\end{equation}
where the indices run through $\{-,+\} = \{-\frac{1}{2},
+\frac{1}{2}\}$. Using Eq.~\eqref{eq:CGSeries1}, the Clebsch-Gordan
decomposition of this representation reads
\begin{equation}
  D_\h^{(\frac{1}{2},\frac{1}{2})} \otimes
  D_\h^{(\frac{1}{2},\frac{1}{2})}
  = D_\h^{(0,0)} \oplus D_\h^{(1,0)} 
  \oplus D_\h^{(0,1)} \oplus D_\h^{(1,1)}\,.
\end{equation}
The subrepresentation by which we have to divide the free algebra $\mathbb{C}\langle
X_{mm'} \rangle[[\h]]$ for the right noncommutative limit as $\h
\rightarrow 0$ is $D_\h^{(1,0)} \oplus D_\h^{(0,1)}$. This corresponds
to the quadratic relations
\begin{equation}
\label{eq:EuclidRels}
\begin{aligned}
  \sum_{m_1,m_2,m_1',m_2'} \CGqs{1/2}{1/2}{1}{m_1}{m_2}{m}
  \CGqs{1/2}{1/2}{0}{m_1'}{m_2'}{0} X_{m_1 m_1'} X_{m_2 m_2'} &= 0 \\
  \sum_{m_1,m_2,m_1',m_2'} \CGqs{1/2}{1/2}{0}{m_1}{m_2}{0}
  \CGqs{1/2}{1/2}{1}{m_1'}{m_2'}{m} X_{m_1 m_1'} X_{m_2 m_2'} &= 0 \,,  
\end{aligned}
\end{equation}
where $m$ runs through $\{-1,0,1\}$. Denoting the generators by
$(\begin{smallmatrix} a & b \\ c & d \end{smallmatrix}) := (X_{mm'})$,
i.e., $d = X_{++}$ etc., relations~\eqref{eq:EuclidRels} read
\begin{equation}
\label{eq:MqRel}
\begin{gathered}
  ab=qba\,, \quad ac=qca\,, \quad bd = q db\,, \quad cd = q dc \\
  bc = cb\,, \quad ad - da = (q - q^{-1}) bc \,,
\end{gathered}
\end{equation}
which are the well known relations of the algebra of 2$\times$2
quantum matrices \cite{Woronowicz:1987b}. The quantum determinant
\begin{equation}
  \mathrm{det}_q := ad -q bc \,,
\end{equation}
is scalar, $(g \otimes g') \tr \mathrm{det}_q = \varepsilon(g \otimes
g') \, \mathrm{det}_q$, and commutes with all generators.

\begin{Definition}
  The $\h$-adic algebra freely generated by $\{a,b,c,d\}$ with
  commutation relations~\eqref{eq:MqRel} is called the $\h$-adic
  quantum Euclidean 4-space or the $\h$-adic algebra of 2$\times$2
  quantum matrices $M_\h(2)$.
\end{Definition} 
Quantum Euclidean 4-space and $M_q(2)$ are the same algebras, for
$\SUh := M_\h(2)/\lrAngle{\mathrm{det}_q=1}$ is Hopf dual to $\Uhsu$
in the sense of \cite{Takeuchi:1977}, which implies that the comodule
algebras of $\SUh$ are the module algebras of $\Uhsu$. In fact, let
$\Delta$ be the coproduct of $M_\h(2)$, $\Delta(X_{ik}) = \sum_j
X_{ij} \otimes X_{jk}$, let $T:M_\h(2) \rightarrow M_\h(2)$ be the
transposition homomorphism which is defined on the generators by
$(X_{ij})^T := X_{ji}$, let $\pi: M_\h(2) \rightarrow
M_\h(2)/\lrAngle{\mathrm{det}_q=1}$ be the canonical epimorphism,
$U^{i}{}_{j} := \pi(X_{ij})$ the generators of $\SUh$, and $\tau$ the
flip of tensor factors. Then the map $\varphi: M_\h(2) \rightarrow
M_\h(2) \otimes \SUh \otimes \SUh$ defined as
\begin{equation}
\label{eq:XComodule}
  \varphi :=  [\id \otimes (\pi \circ T)\otimes \pi ]
    \circ \tau_{12} \circ \Delta^{(2)}
    \quad\Rightarrow\quad \varphi(X_{ij}) = X_{i'j'} \otimes
  U^{i'}{}_{i} \otimes U^{j'}{}_{j}
\end{equation}
is a homomorphism of algebras because it is a concatenation of
homomorphisms, and a corepresentation because $\Delta(U^i{}_k) =
U^i{}_j \otimes U^j{}_k$. Hence, $M_\h(2)$ together with $\varphi$ is
a comodule algebra of $SO_\h(4) := \SUh \otimes \SUh$. The dual of
this coaction~\eqref{eq:XComodule} is the action~\eqref{eq:XModule} by
which $M_q(2)$ becomes a $\Uhso$-module algebra.

A simple ansatz shows, that the only homogeneous highest weight
vectors of $M_\h(2)$ are proportional to $\mathrm{det}_q^k d^l$. In
analogy to the quantum plane, we have to define the minimal degree
irreducible weight vectors by
\begin{equation}
\label{eq:MinkBasis1}
\begin{split}
  T^{(j,j)}_{mm'} 
  &:= q^{\frac{1}{2}[(j-m)(2m-j+1)+(j-m')(2m'-j+1)]} \\
  &\qquad\times
  \sqrt{\frac{[j+m]![j +m']!}{[2j]!^2 [j-m]![j-m']!}}
  \,\, [(F^{j-m} \otimes F^{j-m'}) \tr d^{2j} ] \,,
\end{split}
\end{equation}
for $j \in \tfrac{1}{2}\mathbb{N}_0$, such that they carry a
spin-$(j,j)$ representation,
\begin{equation}
  (g\otimes g') \tr T^{(j,j)}_{mm'}
  = T^{(j,j)}_{\tilde{m}\tilde{m}'}\,
  \rho_\h^{j}(g)^{\tilde{m}}{\!}_{m} \,
  \rho_\h^{j}(g')^{\tilde{m}'}{\!}_{m'} \,.  
\end{equation}
The explicit calculation of Eq.~\eqref{eq:MinkBasis1} leads to
\begin{equation}
\label{eq:MinkBasis2}
\begin{split}
  T^{(j,j)}_{mm'}
  &= \sum_{k} q^{k(m'-m-k)} \qBinom{j-m}{k}_{q^{-2}} 
    \qBinom{j+m'}{j+m'-k}_{q^{-2}}
    \qBinom{2j}{j+m}_{q^{-2}}^{\frac{1}{2}}
    \qBinom{2j}{j+m'}_{q^{-2}}^{-\frac{1}{2}}\\
  &\qquad\times a^{j-m-k} b^k c^{m-m'+k} d^{j+m'-k} \,,
\end{split}
\end{equation}
which reduces $M_\h(2)$ into irreducible subrepresentations by
\begin{equation}
\label{eq:PeterWeyl}
  M_\h(2) = \bigoplus_{j\in \mathbb{N}_0/2} \,
  \bigoplus_{k \in \mathbb{N}_0}\,
  \operatorname{Span}_{\mathbb{C}[[\h]]}
  \{ \mathrm{det}_q^k\, T^{(j,j)}_{mm'} \,|\,  m,m' = -j,\ldots,j\} \,.
\end{equation}
Note that mapping this equation by the canonical epimorphism $\pi$
onto $\SUh$ yields the quantum Peter-Weyl decomposition of $\SUh$ (see
\cite{Schmuedgen}, Sec.~4.2.5). As we argued for the quantum plane,
the reduction of the product of two irreducible weight vectors with
$q$-Clebsch-Gordan coefficients must again be an irreducible weight
vector of the same degree as the product,
\begin{equation}
\label{eq:MinkBasis3}
  \sum_{m_1,m_2, m_1',m_2'}
  \CGqs{j_1}{j_2}{j}{m_1}{m_2}{m}
  \CGqs{j_1'}{j_2'}{j'}{m_1'}{m_2'}{m'}
  T^{(j_1,j_1)}_{m_1m_1'} T^{(j_2,j_2)}_{m_2m_2'}
  = \delta_{jj'} \beta_{j_1 j_2 j} \,
  \mathrm{det}_q^{j_1+j_2-j} T^{(j,j)}_{mm'}
\end{equation}
where the $\beta_{j_1 j_2 j} \in \mathbb{C}[[\h]]$ are $\h$-adic
scalar coefficients. These coefficients can be easily computed by
applying the counit $\varepsilon$ of $M_\h(2)$, for which we have
\begin{equation}
  \varepsilon(T^{(j,j)}_{mm'}) = \delta_{mm'} \,,\qquad
  \varepsilon(\mathrm{det}_q) = 1 \,,
\end{equation}
to Eq.~\eqref{eq:MinkBasis3}. This yields $\beta_{j_1 j_2 j} = 1$ for
all $j_1$, $j_2$, $j$. Finally, moving the $q$-Clebsch-Gordan
coefficients to the other side of Eq.~\eqref{eq:MinkBasis3} produces
the desired expression for the product $\mu_\h$ of $M_\h(2)$,
\begin{equation}
\label{eq:MinkProd1}
  \mu_\h\bigl(
  T^{(j_1,j_1)}_{m_1m_1'} \otimes T^{(j_2,j_2)}_{m_2m_2'} \bigr)
  = \sum_{m, m',j} \CGqs{j_1}{j_2}{j}{m_1}{m_2}{m}
  \CGqs{j_1'}{j_2'}{j}{m_1'}{m_2'}{m'} \,
  \mathrm{det}_q^{j_1+j_2-j} \, T^{(j,j)}_{mm'} \,.
\end{equation}
Again, we compare this with the representations~\eqref{eq:Frep1} and
\eqref{eq:Frep3} of the Drinfeld twists of $\Uhsu$ and obtain
\begin{Proposition}
\label{th:EuclidF}
  Let $\mu_\h$ be the multiplication map~\eqref{eq:MinkProd1} of
  $\h$-adic Euclidean 4-space $M_\h(2)$, $\mu = \lim_{\h \rightarrow
    0} \mu_\h$ its undeformed limit, and $\Fs$ the standard twist of
  Prop.~\ref{th:Frep2}. Then $\mu_\h$ is the
  deformation~\eqref{eq:StarProd2} of $\mu$ by $\F_{\!\so} := \Fs_{13}
  \Fs_{24}$.
\end{Proposition}

Now we can apply Cor.~\ref{th:MinkowskiTwist} in order to obtain the
the multiplication map and twist of quantum Minkowski space, twisting
once more by $\R_{23}^{-1}$. Of course, the multiplication map of
quantum Minkowski space will reproduce the well known commutation
relations of \cite{Carow-Watamura:1990}. For the twist we get
 
\begin{Proposition}
\label{th:MinkowskiF}
Let $\mu_\h$ be the multiplication map of $\h$-adic quantum Minkowski
space $\Mink$, $\mu = \lim_{\h \rightarrow 0} \mu_\h$ its undeformed
limit, and $\Fs$ the standard twist of Prop.~\ref{th:Frep2}. Then
$\mu_\h$ is the deformation~\eqref{eq:StarProd2} of $\mu$ by
$\F_{\!\slC}:= \R^{-1}_{23} \Fs_{13} \Fs_{24}$.
\end{Proposition}

Finally, we consider real structures. Since $\Fs$ is real in the sense
of Eq.~\eqref{eq:Freal}, so is $\F_{\!\so}$. Moreover, since
$\F_{\!\slC}$ is the twist of $\F_{\!\so}$ by $\R_{23}^{-1}$, which
was shown to be real in the proof of Cor.~\ref{th:MinkowskiTwist},
(iv), $\F_{\!\slC}$ is real, as well. Hence,
Prop.~\ref{th:StarModuleDeform} applies to both, quantum Euclidean
4-space and quantum Minkowski space. From a reasoning which is
completely analogous to the one that led to Eq.~\eqref{eq:PlaneStar}
we conclude that the deformation $*_\h = * + \Ocal(\h)$ of the real
structure $*$ of the undeformed Minkowski spacetime algebra $\Xcal
\equiv \Mink$ has to be defined by
\begin{equation}
  x^{*_\h} =
  (\E^{-\frac{\h H}{2}} \otimes \E^{-\frac{\h H}{2}}) \tr x^*
\end{equation}
such that $(\Xcal, \mu_\h, *_\h)$ becomes a module $*$-algebra of
$\UslC[[\h]]$ with respect to the deformed Hopf structure.

\section{Conclusion}

It is possible to use Drinfeld twists in order to realize quantum
spaces as covariant star products. We have shown this for three
important examples, the quantum plane (Prop.~\ref{th:PlaneF}), quantum
Euclidean 4-space (Prop.~\ref{th:EuclidF}), and quantum Minkowski
space (Prop.~\ref{th:MinkowskiF}). While it was known that the
Drinfeld twists control the deformation of enveloping algebras into
quantum enveloping algebras, it is now clear that certain twists also
control the deformation of spaces into quantum spaces. This is not
unexpected, since the covariance condition of the action of a symmetry
on a space algebra ties the Hopf structure of the symmetry algebra
closely to the multiplicative structure of the space algebra. Our
considerations included real structures of quantum enveloping algebras
and quantum spaces. In Prop.~\ref{th:StarModuleDeform} we have
formulated a sufficient condition on the Drinfeld twist, its reality
in the sense of Eq.~\eqref{eq:Freal}, to be compatible with the real
structure of a quantum space, and we have shown that there is a unique
element $\sigma$ of the enveloping algebra which implements the
deformation of the real structure.

Star products are often defined by identifying the vector spaces of
two space algebras $\Xcal$ and $\Xcal_\h$ by an vector space
isomorphism, $\varphi: \Xcal \rightarrow \Xcal_\h$, and transferring
the multiplication by $x \star y :=
\varphi^{-1}[\varphi(x)\varphi(y)]$. The Moyal-Weyl product is an
example of this procedure. The linear isomorphism $\varphi$ is called
an ordering prescription, because it defines how an ordered monomial
of the commutative algebra $\Xcal$ has to be represented in the
noncommutative algebra $\Xcal_\h$. For example $\varphi(xy) =
\frac{1}{2}[\hat{x}\hat{y} + \hat{y}\hat{x}]$, where $\hat{x} :=
\varphi(x)$, $\hat{y} := \varphi(y)$, for the symmetric ordering. The
star product which is obtained by the standard twist of
Prop.~\ref{th:PlaneF} amounts to the ordering prescription which
identifies the basis vectors which completely reduce the space and the
quantum space respectively into its irreducible subrepresentations,
which is a natural ordering in the context of representation theory.
For the quantum plane this is almost the lexicographic
(Poincar\'{e}-Birkhoff-Witt) ordering,
\begin{equation}
  \varphi(x^k y^l) =
  \qBinom{k+l}{k}_{q^{-2}}^{\frac{1}{2}}
  \tbinom{k+l}{k}^{-\frac{1}{2}}
  x^k y^l \,,
\end{equation}
where we recall that $x \equiv x_-$, $y \equiv x_+$. The basis of the
spin-$j$ subrepresentation is unique up to scalar multiples. A
rescaling $T^j_m \mapsto \beta(j) T^j_m$ with $\beta(j) = 1 +
\Ocal(\h)$ would change the multiplication map by multiplying the
right hand side of Eq.~\eqref{eq:qPlaneMulti} with $\beta(j)
\beta^{-1}(j_1) \beta^{-1}(j_2)$. Identifying the scale factors with
the representations $\beta(j) = \rho^j(z)$ of some invertible central
element $z$, the twist which realizes the star product must be
redefined by $\F \mapsto \F (z \otimes z) \Delta z^{-1}$. We conclude
that the twist of Prop.~\ref{th:PlaneF} which realizes the star
product of the quantum plane is unique up to a central 2-coboundary.
The standard twist is the unique twist with the additional property
$y^{\star n} = y^n$. An analogous statement is true for quantum
Euclidean 4-space and quantum Minkowski space.

To our knowledge, no Drinfeld twist for the Drinfeld-Jimbo deformation
of a semisimple Lie algebra had so far been computed explicitly (see
\cite{Fiore:1997} for the Heisenberg algebra). We have circumvented
this problem by reducing the algebraic questions to representation
theoretic ones. Although $\Uhsu$ is the simplest case conceivable, an
algebraic order by order calculation of the twist runs quickly into
overwhelming combinatorial problems \cite{Dabrowski:1996}. An
alternative approach would be the reconstruction of the twist from its
representations, which would profit from the computational effort that
has gone into the calculation of the $q$-Clebsch-Gordan coefficients.


\appendix

\section{Appendix}

\textit{Proof of Corollary~\ref{th:HopfTwist}.} Let throughout the proof
$g \in \Ug[[\h]]$ be an arbitrary element of the $\h$-adic enveloping
algebra.
  
(i) $\F \Delta(g) \F^{-1} = \Delta_\h(g) = \F' \Delta(g) \F^{\prime
  -1}$ implies that $\F^{-1}\F'$ commutes with $\Delta(g)$ for all
$g$. Conversely, let $\Tcal \Delta(g) \Tcal^{-1} = \Delta(g)$. Then
$\F \Tcal \Delta(g) \Tcal^{-1} \F^{-1} = \F \Delta(g) \F^{-1} =
\Delta_\h(g)$ for all $g$.
  
(ii) By the left counit property of $\varepsilon_\h$ we get
\begin{multline}
  \varepsilon(g) = \varepsilon \bigl( \varepsilon_\h(g_{(1_\h)})
  g_{(2_\h)} \bigr) = \varepsilon\bigl( \varepsilon_\h
  (\F_{[1]})\F_{[2]} \, \varepsilon_\h(g_{(1)}) g_{(2)} \,
  \varepsilon_\h (\F^{-1}_{[1']})\F^{-1}_{[2']} \bigr) \\
  = \varepsilon_\h(g_{(1)}) \varepsilon(g_{(2)}) \, \varepsilon_\h
  (\F_{[1]})\varepsilon(\F_{[2]}) \varepsilon_\h
  (\F^{-1}_{[1']})\varepsilon(\F^{-1}_{[2']})
      = \varepsilon_\h\bigl( g_{(1)} \varepsilon(g_{(2)})\bigr) 
      = \varepsilon_\h(g).
\end{multline}
  
(iii) Let us define the left and right counit constraints by $l :=
\varepsilon(\F_{[1]})\F_{[2]}$ and $r := \F_{[1]}
\varepsilon(\F_{[2]})$. From the left counit property of
$\varepsilon_\h = \varepsilon$ it follows that $g = \varepsilon
(g_{(1_\h)}) g_{(2_\h)} = \varepsilon(\F_{[1]})\F_{[2]} \,
\varepsilon(g_{(1)}) g_{(2)} \,
\varepsilon(\F^{-1}_{[1']})\F^{-1}_{[2']} = l g l^{-1}$ and
analogously for the right counit property $g = rgr^{-1}$. Hence,
$\Tcal := \varepsilon(l)(r^{-1} \otimes l^{-1})$ is
$\mathfrak{g}$-invariant. By (i) $\F' := \F\Tcal$ is a twist with $l'
:= \varepsilon(\F'_{[1]})\F'_{[2]} = \varepsilon(l)
\varepsilon(\F_{[1]}r^{-1}) \F_{[2]} l^{-1} = \varepsilon(l r^{-1}) l
l^{-1} = 1$, where we have used $\varepsilon(l) =
\varepsilon(\F_{[1]}\F_{[2]}) = \varepsilon(r)$. Analogously, we find
that $r' = 1$. Since $\F'$ is invertible, $\F' = \beta 1 + \Ocal(\h)$
for some complex number $\beta \neq 0$, and $1= l' = \beta +
\Ocal(\h)$ it follows that $\beta = 1$.
  
(iv) Define $\sigma_1 := S_\h(\F_{[1]})\F_{[2]}$. Then
\begin{multline}
  S_\h(g) \sigma_1 = S_\h(g_{(1)}) \sigma_1 g_{(2)} S(g_{(3)})
  = S_\h(\F_{[1]} g_{(1)(1)}) \F_{[2]} g_{(1)(2)} S(g_{(2)})\\
  = S_\h(g_{(1)(1_\h)}\F_{[1]}) g_{(1)(2_\h)} \F_{[2]} S(g_{(2)}) =
  S_\h(\F_{[1]}) \varepsilon(g_{(1)}) \F_{[2]} S(g_{(2)}) = \sigma_1
  S(g)\,,
\end{multline}
where we have used the left coinverse property of $S_\h$. Analogously,
defining $\sigma_2^{-1} := \F_{[1]}^{-1} S_\h(\F_{[2]}^{-1})$ and
using the right coinverse property of $S_\h$ we get $\sigma_2^{-1}
S_\h(g) = S(g) \sigma_2^{-1}$. Since $\F$ is invertible, so are
$\sigma_1$ and $\sigma_2^{-1}$. Thus, $S_\h$ and $S$ are related by
the inner automorphisms $\sigma_1 S(g) \sigma_1^{-1} = S_\h(g) =
\sigma_2 S(g) \sigma_2^{-1}$. The antipode is surjective, so
$\sigma_1^{-1} \sigma_2$ must be central. Moreover, $(\sigma_1^{-1}
\sigma_2) (\sigma_1 \sigma_2^{-1}) = \sigma_1 (\sigma_1^{-1} \sigma_2)
\sigma_2^{-1} = 1$, hence, $\sigma_1^{-1}\sigma_2 =
\sigma_2\sigma_1^{-1}$. Finally,
\begin{equation}
    \sigma_1^{-1} = \sigma_1^{-1}
    S_\h(\F_{[1']}\F_{[1]}^{-1}) \F_{[2']} \F_{[2]}^{-1} 
    = \sigma_1^{-1} S_\h(\F_{[1]}^{-1}) \sigma_1 \F_{[2]}^{-1} 
    = S(\F^{-1}_{[1]}) \F^{-1}_{[2]} \,,
\end{equation}
and, analogously, $\sigma_2 = \F_{[1]} S(\F_{[2]})$.
 
(v) From the coassociativity of $\Delta_\h$ we deduce
\begin{multline}
    \bigl[(\Delta_\h \otimes \id) \circ \Delta_\h \bigr](g) 
    = \F_{12}(\Delta \otimes \id)(\F)\, (\Delta^{\!(2)}g) \,
    (\Delta \otimes \id)(\F^{-1}) \F^{-1}_{12}\\
    = \bigl[(\id \otimes \Delta_\h) \circ \Delta_\h\bigr](g)
    = \F_{23} (\id \otimes \Delta )(\F) \, (\Delta^{\!(2)}g) \,
    (\id \otimes \Delta)(\F^{-1}) \F^{-1}_{23} \,.
\end{multline}
Hence, the coassociator~\eqref{eq:CoassDef} commutes with all
$\Delta^{\!(2)}g$.
  
(vi) Since every bialgebra isomorphism is automatically a Hopf algebra
isomorphism and since by (ii) the counit stays undeformed, we only
have to consider the coproducts.  Let $\Delta$ and $\Delta_\h$ be
related by a twist $\F = (u \otimes u)\Delta u^{-1}$. Then $g \mapsto
u g u^{-1}$ maps $\Delta$ to $\Delta_\h$. Conversely, assume that the
two Hopf structures are isomorphic, so there is an algebra
automorphism $\beta$ with $\Delta_\h = (\beta \otimes \beta) \circ
\Delta \circ \beta^{-1}$. Since $\Delta_\h = \Delta + \Ocal(\h)$ we
can choose $\beta = \id + \Ocal(\h)$. Since $\mathfrak{g}$ is
semisimple we have $H^1(\Ug,\Ug) = 0$, which implies that this
automorphism is inner, $\beta(g) = ugu^{-1}$. \qed

\end{document}